\documentclass{amsart}
\usepackage{amsmath, amssymb}

\let\a\alpha
\let\b\beta
\let\g\gamma
\let\f\varphi
\let\k\kappa
\let\l\lambda
\let\L\Lambda
\let\m\mu

\let\x\xi
\let\z\zeta
\let\o\omega

\def\E{\mathcal E}
\def\F{\mathcal F}

\def\Gtilda{\widetilde{G}}

\def\Ltilda{\widetilde{\Lambda}}
\def\O{\mathcal O}
\def\P{\mathbb P}
\def\Wtilda{\widetilde{W}}

\newcommand{\lra}{\longrightarrow}
\newcommand{\tensor}{\otimes}
\newcommand{\isom}{\simeq}
\newcommand{\bdm}{\begin{displaymath}}
\newcommand{\edm}{\end{displaymath}}
\newcommand{\Aut}{\text{Aut}}
\newcommand{\Hom}{\text{Hom}}
\newcommand{\GL}{\text{GL}}
\newcommand{\Ker}{\text{Ker}}

\begin{document}

\title{On some Quotients Modulo Nonreductive Groups}

\author{Mario Maican}

\maketitle

\tableofcontents

\section{Introduction}

Quotients modulo nonreductive groups arise in the study of moduli spaces
of sheaves. Let $k$ be an algebraically closed field of characteristic zero.
We fix coherent algebraic sheaves $\E$ and $\F$ on the projective space
$\P^r=\P^r(k)$, which are direct sums of simple sheaves. The linear algebraic
group $G=\Aut(\E) \times \Aut(F)$ acts by conjugation on the finite dimensional
vector space $W= \text{Hom}(\E,\F)$. To every character of $G$ Dr\'ezet and
Trautmann associate subsets $W^{ss} \subset W$ of \emph{semistable points}
and $W^s \subset W$ of \emph{stable points}. One expects there to be
categorical quotients $W^{ss}//G$ and $W^s/G$. In \cite{maican} the author
has found many examples of locally closed subsets inside moduli spaces
of semistable sheaves (in the sense of Gieseker) on $\P^2$, with linear Hilbert
polynomial, that are isomorphic to quotients of the form $W^{ss}//G$.

This paper is not concerned with the study of moduli spaces but only with
the question of existence of quotients modulo $G$. If $\E$ has more than one
kind of simple sheaf in its decomposition, then $\text{Aut}(\E)$ is
nonreductive, so this situation falls outside Mumford's Geometric Invariant
Theory. Dr\'ezet and Trautmann address this difficulty in \cite{drezet-trautmann}
by embedding the action of $G$ on $W$ into the action of a reductive group
$\Gtilda$ onto a larger space $\Wtilda$. Their main result states that if certain
compatibility conditions relating the semistable points of $W$ and $\Wtilda$
are satisfied, then $W^{ss}//G$ and $W^s/G$ exist and are projective,
respectively quasiprojective varieties. Moreover, they find sufficient conditions,
expressed in terms of linear algebra constants, under which the compatibility
conditions are fulfilled. This allows Dr\'ezet and Trautmann to establish the
existence of quotients for certain classes of morphisms, notably for
morphisms of the form
\bdm
m_1 \O(-2) \oplus m_2 \O(-1) \lra n\O,
\edm
cf. 6.4 in \cite{drezet}. The purpose of this paper is to give more
examples to the Dr\'ezet-Trautmann theory. We use their embedding into the
action of $\Gtilda$ but we do not use their linear algebra constants. We are
concerned only with the geometric quotients $W^s/G$ and we do not discuss
properly semistable morphisms, i.e. morphisms which are semistable but not
stable. Applying 6.6.1 from \cite{drezet-trautmann} we establish the existence
of geometric quotients in the following situations:
\begin{enumerate}
\item[$\bullet$] $m\O(-d_1) \oplus 2\O(-d_2) \lra n\O$ on $\P^r$ such that
$0 < \l_1 < 1/(2a+m)$ and either the conditions
\bdm
m < {r-1+d_1-d_2 \choose r-1}, \quad \l_1 \le \frac{1}{a+m-1} \left( 1 - \frac{1}{n}
{r+d_2 -1 \choose r} \right),
\edm
or the conditions
\bdm
m < {r+d_1-d_2 \choose r}, \quad \l_1 \le \frac{1}{3m} - \frac{2}{3mn}
{r+d_2-1 \choose r}
\edm
are satisfied. Here $m$ and $n$ are not both even;
\item[$\bullet$] $\O(-d-1) \oplus 3\O(-d) \lra n\O$ on $\P^2$ such that
\begin{align*}
0 & < \l_1 < \frac{1}{10}, \\
-\frac{1}{2} + \frac{3}{4n}(d^2+d) & \le \l_1 \le \frac{2}{5} - \frac{3}{10 n}
(d^2 + d), \\
-2 + \frac{3}{n}(d^2+2d) & \le \l_1 \le 1- \frac{3}{4n}(d^2+3d);
\end{align*}
\item[$\bullet$] $m\O(-d-1) \oplus 3\O(-1) \lra n\O$ on $\P^2$ with $m < a$ and
\bdm
0 < \l_1 < \frac{1}{3a+m}, \quad \l_1(4m-3a+3b) \le \frac{n-3}{n}, \quad
\l_1 \le \frac{n-6}{mn}.
\edm
Here $a=(d+1)(d+2)/2$ and $b=d(d+1)/2$ and $3 \nmid \text{gcd}(m,n)$;
\item[$\bullet$] $\O(-d-2)\oplus 3\O(-d) \lra n\O$ on $\P^2$ such that
\begin{align*}
\frac{6}{19} < & \l_2 < \frac{1}{3}, \\
& \l_2 \le \frac{1}{2} - \frac{1}{2n} {d+1 \choose 2}, \\
& \l_2 \le 1-\frac{1}{n} {d+2 \choose 2} - \frac{1}{n} {d+1 \choose 2}, \\
& \l_2 \le \frac{1}{n} {d+2 \choose 2};
\end{align*}
\item[$\bullet$] $m\O(-d_3) \oplus \O(-d_2) \oplus \O(-d_1) \lra n\O$ on $\P^r$
satisfying
\begin{align*}
m & <{d_1-d_2 +r-1 \choose r-1}, \\
a_{21} \l_1 <  \l_2 & < \frac{1-m\l_1-a_{31}\l_1+a_{32}a_{21}\l_1}{1+a_{32}}, \\
\l_2 & < 1-m\l_1, \\
\l_1 & \le \frac{1}{m+a_{21}} - \frac{1}{mn+ a_{21}n} {d_3 + r \choose r}.
\end{align*}
\end{enumerate}
The precise meaning of $\l_1$ and $\l_2$ is revealed in section 2: they encode
the character of $G$ relative to which the sets of stable points have been defined.
Beside the above, we have at (3.3) a more general criterion giving sufficient
conditions, expressed in terms of linear algebra constants, under which
morphisms of type (2,1)
\bdm
m_1 \O(-d_1) \oplus m_2 \O(-d_2) \lra n\O
\edm
satisfy the compatibility conditions from \cite{drezet-trautmann} that lead to
existence of geometric quotients.

The paper is organized as follows: in section 2 we supply background material
about quotients modulo nonreductive groups. As this is not part of mainstram
Geometric Invariant Theory, we felt it necessary to reproduce the main results
and definitions from \cite{drezet-trautmann}. Section 3 contains our criterion
for morphisms of type (2,1). We apply this criterion in section 7. The remaining
sections are indirect applications of (3.3), by which we mean applications of the
method of proof, rather than the statement of (3.3). All polarizations we discuss
are assumed to be nonsingular in a sense that is made precise in our very
brief discussion of the Geometric Invariant Theory fan, cf. the beginning of
section 3.


\section{Quotients Modulo Nonreductive Groups}

In this section we reproduce some notations and definitions and we quote
the main result from \cite{drezet-trautmann}. We also quote two particular cases
of King's Criterion of Semistability, as formulated in \cite{drezet-trautmann},
which we will use in the subsequent sections.

Fix a vector space $V$ over $k$ of dimension $n+1$. Dr\'ezet and Trautmann consider
coherent algebraic sheaves $\E$ and $\F$ on the projective space $\P^n=\P(V)$
having decompositions
\bdm
\E = \bigoplus_{1 \le i \le r} M_i \tensor \E_i, \qquad
\F = \bigoplus_{1 \le l \le s} N_l \tensor \F_l.
\edm
Here $M_i$, $N_l$ are vector spaces over $k$ of dimensions $m_i$, $n_l$.
In \cite{drezet-trautmann} it is assumed that $\E_i$ and $\F_l$ are simple sheaves,
but for the purposes of this paper we will assume that they are line bundles:
\bdm
\E_i = \O(e_i), \quad e_1 < \ldots < e_r, \qquad \F_l = \O(f_l), \quad f_1 < \ldots < f_s.
\edm
The linear algebraic group $\Aut(\E) \times \Aut(\F)$ acts by conjugation on the finite
dimensional vector space
\bdm
W = \text{Hom}(\E,\F).
\edm
The subgroup of homotheties, which we identify with $k^*$, acts trivially so, without
losing any information, we can instead consider the action of the quotient
\bdm
G= \Aut(\E) \times \Aut(\F)/k^*.
\edm
If $r>1$, or if $s>1$, the group $G$ is nonreductive, however, it contains the reductive subgroup
\bdm
G_{\text{red}}= \GL (M_1) \times \ldots \times \GL(M_r) \times \GL(N_1) \times \ldots
\times \GL(N_s)/k^*.
\edm
We represent elements of $G_{\text{red}}$ by pairs $(g,h)$, with
\bdm
g=(g_1,\ldots,g_r), \quad h=(h_1,\ldots,h_s), \quad g_i \in \GL(M_i), \quad h_j \in \GL(N_j).
\edm
The characters of $G_{\text{red}}$ are of the form
\bdm
\chi(g,h) = \prod_{1 \le i \le r} \text{det}(g_i)^{-\l_i} \cdot \prod_{1 \le l \le s} \text{det}(h_l)^{\m_l}
\edm
for integers $\l_1,\ldots,\l_r,\m_1,\ldots,\m_s$. As $\chi$ must be trivial on the subgroup
of homotheties, we require the condition
\bdm
\sum_{1 \le i \le r} m_i \l_i = d = \sum_{1 \le l \le s} n_l \m_l.
\edm
Clearly $\chi$ extends to a character of $G$, which we also denote by $\chi$.
Dr\'ezet and Trautmann call a \emph{polarization} the tuple
\bdm
\L = (-\l_1/d,\ldots,-\l_r/d,\m_1/d,\ldots,\m_s/d).
\edm
They consider a semistability notion for $G_{\text{red}}$ depending on $\L$.
We quote below the definition. We point out that, relative to a certain linearization,
this is the usual notion of semistability from Geometric Invariant Theory.
This is made precise at lemma 3.4.1 in \cite{drezet-trautmann}.

\noindent \\
{\bf (2.1) Definition:} Let $\L$ be a fixed polarization. A point $\f \in W$ is called
\begin{enumerate}
\item[(i)] \emph{semistable} with respect to $G_{\text{red}}$ and $\L$ if there are
an integer $m \ge 1$ and an algebraic function $f$ on $W$ satisfying
$f(g.w)=\chi^m(g)f(w)$ for all $g \in G_{\text{red}}$ and $w \in W$, such that $f(\f) \neq 0$;
\item[(ii)] \emph{stable} with respect to $G_{\text{red}}$ and $\L$ if the isotropy group
of $\f$ in $G_{\text{red}}$ is finite and there is $f$ as above, but with the additional
property that the action of $G_{\text{red}}$ on the set $\{ w \in W,\ f(w) \neq 0 \}$
has closed orbits.
\end{enumerate}

\noindent \\
This definition is consistent because proportional tuples of integers give rise to the
same sets of semistable (stable) points. Let $T$ be the maximal torus of $G_{\text{red}}$.
Note that $T$ is also a maximal torus in $G$. A point $\f \in W$ is semistable (stable)
with respect to $G_{\text{red}}$ and $\L$ if and only if every point in its
$G_{\text{red}}$-orbit is semistable (stable) with respect to $T$ and the restriction
of $\chi$ to $T$. Taking this equivalence as definition for semistability (stability)
with respect to $G$ we arrive at the following concept introduced by Dr\'ezet 
and Trautmann:

\noindent \\
{\bf (2.2) Definition:} A point $\f \in W$ is called \emph{semistable} (\emph{stable})
with respect to $G$ and $\L$ if every point in its $G$-orbit is semistable (stable)
with respect to $G_{\text{red}}$ and $\L$. We denote by $W^{ss}(G,\L)$ and
$W^s(G,\L)$ the sets of semistable, respectively stable points in $W$. \\

For checking semistability in concrete situations we need a criterion derived by
A. King in \cite{king} from Mumford's Numerical Criterion. We use its formulation from
\cite{drezet-trautmann}. Below we quote only a particular case. Let us represent
a point $\f \in W$ by a matrix
\bdm
(\f_{li})_{1 \le l \le s, 1 \le i \le r} \quad \text{with} \quad \f_{li} \in \text{Hom}
(M_i \tensor H_{li}^*,N_l), \quad H_{li}= \text{Hom}(\E_i,\F_l).
\edm
A family of subspaces $M_i' \subset M_i$, $N_l' \subset N_l$ will be called
\emph{admissible} if not all subspaces are zero and we do not have
$M_i'=M_i$, $N_l'=N_l$ for all $i$, $l$.

\noindent \\
{\bf (2.3) Proposition:} \emph{A morphism $\f \in W$ is semistable (stable) with respect
to $G$ and $\L$ if and only if for each admissible family of subspaces $M_i' \subset M_i$,
$N_l' \subset N_l$, which satisfies}
\bdm
\f_{li} (M_i' \tensor H_{li}^*) \subset N_l' \qquad \text{\emph{for all}} \quad i,\, l,
\edm
\emph{we have}
\bdm
\sum_{l=1}^s \m_l \, \text{\emph{dim}}(N_l') \, \ge \, (>) \, \sum_{i=1}^r \l_i \, \text{\emph{dim}}(M_i'). \\
\edm

Dr\'ezet and Trautmann embed the action of $G$ on $W$ into the action of a reductive
group $\Gtilda$ on a finite dimensional vector space $\Wtilda$. They introduce associated
polarizations
\bdm
\Ltilda = (\a_1,\ldots,\a_r,\b_1,\ldots,\b_s)
\edm
and define sets of semistable and stable points $\Wtilda^{ss}(\Gtilda,\Ltilda)$,
respectively $\Wtilda^{s}(\Gtilda,\Ltilda)$, as at (2.1). In general they determine
the following relationships:

\noindent \\
{\bf (2.4) Proposition:} \emph{Let $\z: W \lra \Wtilda$ denote the embedding map.
Then we have the inclusions}
\bdm
\z^{-1} (\Wtilda^{ss}(\Gtilda,\Ltilda)) \subset W^{ss}(G,\L)
\qquad \text{\emph{and}} \qquad
\z^{-1} (\Wtilda^{s}(\Gtilda,\Ltilda)) \subset W^{s}(G,\L).
\edm

\noindent \\
They also find sufficient conditions under which the reverse inclusions hold but,
in this paper, we will not use them. The main result from \cite{drezet-trautmann}
states that, if the sets of semistable (stable) points in $W$ and $\Wtilda$ are compatible,
then there are good or geometric quotients modulo $G$. We quote below the part
that will be used in the sequel:

\noindent \\
{\bf (2.5) Proposition:} \emph{If $\z^{-1} (\Wtilda^{s}(\Gtilda,\Ltilda)) = W^{s}(G,\L)$,
then there exists a geometric quotient $W^s(G,\L)/G$ which is a smooth, quasiprojective
variety.}

\noindent \\
Combining (2.4) and (2.5) one arrives at the following:

\noindent \\
{\bf (2.6) Proposition:} \emph{If $\z (W^s(G,\L)) \subset \Wtilda^{s}(\Gtilda,\Ltilda)$,
then there exists a geometric quotient $W^s(G,\L)/G$ which is a smooth, quasiprojective
variety.}

\noindent \\
The goal of this paper is to give classes of examples of geometric quotients that result
as applications of (2.6). We will avoid any discussion of properly semistable morphisms,
i.e. morphisms which are semistable but not stable. Dr\'ezet and Trautmann's theory
works under certain a priori restrictions on the polarization $\L$.
First, it is noticed that the set of semistable (stable) points in $W$ is nonempty only if
\bdm
\l_i \, \ge \, (>) \, 0 \quad \text{and} \quad \m_l \, \ge \, (>) \, 0 \quad \text{for all} \quad
i,\, l.
\edm
Other conditions can be found at 5.4.1 in \cite{drezet-trautmann}: $\Wtilda^{s}(\Gtilda,\Ltilda)$
is not empty only if
\begin{align*}
\tag{2.7} \a_2 & > 0 && \text{for} \quad (r,s)=(2,1), \\
\a_3 & > 0, \, \l_1 p_1 < 1 && \text{for} \quad (r,s)=(3,1).
\end{align*}

\noindent \\
We will describe the embedding only in the case $(r,s)=(2,1)$. Let us write
\bdm
\E_1 = \O(-d_1), \quad \E_2 = \O(-d_2), \quad \F_1 = \O.
\edm
The polarization $\L$ is a triple
\bdm
\L=(-\l_1,-\l_2,\m_1), \quad \text{with} \quad \m_1 = \frac{1}{n} \quad \text{and}
\quad m_1 \l_1 + m_2 \l_2 = 1.
\edm
Recall that
\bdm
H_{11}= \Hom(\E_1,\F_1)= S^{d_1} V^*, \qquad H_{12}= \Hom(\E_2,\F_1)=S^{d_2}V^*.
\edm
Consider also the spaces
\begin{align*}
A_{21} & = \Hom(\E_1,\E_2)= S^{d_1-d_2}V^*, \\
P_1 & = \Hom(\E_1,\E)= M_1 \, \oplus \, M_2 \tensor A_{21}, \\
P_2 & = \Hom(\E_2,\E)= M_2, \\
\Wtilda & = \Hom(P_2 \tensor A_{21}, P_1) \, \oplus \, \Hom(P_1,H_{11} \tensor N_1).
\end{align*}
Small case letters $a=a_{21}$, $p_1=m_1+m_2 a_{21}$, $p_2=m_2$ denote the dimensions
of the corresponding spaces. The group $\Gtilda$ acting by conjugation on $\Wtilda$ is
\bdm
\Gtilda = \GL(P_1) \times \GL(P_2) \times \GL(N_1)/k^*,
\edm
and the associated polarization is
\bdm
\Ltilda = (-\a_1,-\a_2,\b_1) \quad \text{with} \quad \a_1 = \l_1, \ \a_2 = \l_2-a\l_1,\ \b_1 = \m_1 = \frac{1}{n_1}.
\edm
We represent morphisms in $W$ as matrices $\f=(\f',\f'')$,
\begin{align*}
\f' & = (\f_{ij})_{1 \le i \le n_1,\, 1 \le j \le m_1}, && \f_{ij} \in H_{11}, \\
\f'' & = (\f_{ij})_{1 \le i \le n_1,\, m_1 < j \le m_1 + m_2}, && \f_{ij} \in H_{12}.
\end{align*}
We write $\f_j$ for the $j^{\text{th}}$ column of $\f$. We consider the row vector
\bdm
X = \left[
\begin{array}{ccc}
X_1 & \cdots & X_a
\end{array}
\right]
\edm
with entries forming a basis of $A_{21}$. The $p_1 \times p_2$-matrix with entries in $A_{21}$
\bdm
\xi = \left[
\begin{array}{cccc}
0 & 0 & \cdots & 0 \\
X^T & 0 & \cdots & 0 \\
0 & X^T & \cdots & 0 \\
\vdots & \vdots & \ddots & \vdots \\
0 & 0 & \cdots & X^T
\end{array}
\right]
\edm
represents an element in $\Hom(P_2 \tensor A_{21}, P_1)$. The $n_1 \times p_1$-matrix
with entries in $H_{11}$
\bdm
\g(\f) = \left[
\begin{array}{cccc}
\f' & \f_{m_1+1} X & \cdots & \f_{m_1+m_2} X
\end{array}
\right]
\edm
represents an element in $\Hom(P_1,H_{11} \tensor N_1)$. We put $\z(\f) = (\xi, \g(\f))$.
It is clear that $\g (h_1 \f) = h_1 \g(\f)$ for all $h_1 \in \GL(N_1)$. For this reason, when it
comes to semistability considerations, we can, and we will assume that $h_1$ is the identity
automorphism.

We finish this section with the particular case of King's Criterion of Semistability that applies
to $\Wtilda$:

\noindent \\
{\bf (2.8) Proposition:} \emph{A point $(x,\g)$ of $\Wtilda$ is semistable (stable) with respect
to $\Gtilda$ and $\Ltilda$ if and only if for each admissible family of subspaces $P_1' \subset
P_1$, $P_2' \subset P_2$, $N_1' \subset N_1$ satisfying}
\bdm
x(P_2' \tensor A_{21}) \subset P_1', \qquad \g(P_1' \tensor H_{11}^*) \subset N_1',
\edm
\emph{we have}
\bdm
\b_1 \, \text{\emph{dim}}(N_1') \, \ge \, (>) \, \a_1 \, \text{\emph{dim}}(P_1') +
\a_2 \, \text{\emph{dim}}(P_2').
\edm


\section{A Criterion for Morphisms of Type (2,1)}

We fix integers $d_1 > d_2 > 0$, we fix a vector space $V$ of dimension $r+1$ and we
consider morphisms
\bdm
\f = (\f',\f''): m_1 \O(-d_1) \oplus m_2 \O(-d_2) \lra n\O \qquad \text{on} \quad \P^r = \P(V).
\edm
The polarization $\L$ is uniquely determined by $\l_1 \in [0,1/m_1]$. The theory of the
GIT-fan, as developed in \cite{ressayre} and other works, informs us that there are finitely
many values $0 = s_0 < s_1 < \ldots < s_q = 1/m_1$ such that when $\l_1$ varies in an
interval $(s_{\k},s_{\k+1})$ the set of semistable morphisms does not change
and each open interval $(s_{\k},s_{\k+1})$ is maximal with this property.
The intervals $(s_{\k}, s_{\k+1})$ are called \emph{chambers}. The points $s_{\k}$
are called \emph{singular values} for $\l_1$. Either $\l_1$ or $\l_2$ uniquely determine $\L$,
so we can talk of \emph{singular polarizations} $\L$, or \emph{singular values} for $\l_2$.

According to King's Criterion of Semistability (2.3), whenever the set of properly semistable morphisms
with respect to $\L$ is nonempty, there is an equality
\bdm
\sum_{l=1}^s \m_l \, \text{dim}(N_l') = \sum_{i=1}^r \l_i \, \text{dim}(M_i').
\edm
Those polarizations for which there is an equality as above, for some choice of subspaces
$N_l'$ and $M_i'$, will be called \emph{irregular}, or we may say that $\l_1$ or $\l_2$ is
\emph{irregular}. There are situations in which all polarizations are irregular, but those
situations will not be addressed in this paper. The other possibility for morphisms of type (2,1),
which we assume henceforth, is a finite set of irregular polarizations.
From King's Criterion of Semistability we see that $W^s (G,\L)$ depends only on the set
of tuples of integers $(a_1, \ldots, a_r, b_1, \ldots, b_s)$, $0 \le a_i \le m_i$, $0 \le b_l \le n_l$,
for which there is an inequality
\bdm
\sum_{l=1}^s \m_l  b_l > \sum_{i=1}^r \l_i  a_i.
\edm
Let now $\L^0$ be a fixed regular polarization. If $\L$ is sufficiently close to $\L^0$, meaning that
$\l_1$ is sufficiently close to $\l_1^0$, then $\L$ is also regular and the sets of tuples of integers  
with the above property for $\L$ and $\L^0$ are the same. Thus, the sets of stable morphisms,
which for regular polarizations coincide with the sets of semistable morphisms, are the same.
In other words, any open interval bounded by two consecutive irregular values for $\l_1$
is contained in a chamber. Thus, all singular polarizations are irregular.
The author does not know if the converse statement is also true, but he will give below an
example in which the irregular polarizations are singular.
These considerations also show that for $\L$ in a chamber we have $W^{ss}(G,\L)=W^s(G,\L)$ because,
if $\L$ happens to be irregular, we can perturb it slightly to a regular polarization in the same
chamber.

Given integers $0 \le \k_1 \le m_1$ and $0 \le \k_2 \le m_2$ we denote by $l_{\k_1\k_2}$
the smallest integer satisfying
\bdm
\frac{l_{\k_1\k_2}}{n} > \k_1 \l_1 + \k_2 \l_2
\edm
and we consider morphisms of the form $\f_{\k_1\k_2}=(\f_{\k_1 \k_2}',\f_{\k_1 \k_2}'')$, where
\bdm
\f_{\k_1 \k_2}' = \left[
\begin{array}{ll}
\star & 0_{l_{\k_1\k_2},m_1-\k_1} \\
\star & \star
\end{array}
\right], \qquad \f_{\k_1 \k_2}'' = \left[
\begin{array}{ll}
\star & 0_{l_{\k_1\k_2},m_2-\k_2} \\
\star & \star
\end{array}
\right].
\edm
By $0_{lk}$ we denote the identically zero $l \times k$-matrix.
According to King's Criterion (2.3), the morphism $\f$ is semistable if and only if
it is not equivalent to $\f_{\k_1 \k_2}$ for any choice of $\k_1$ and $\k_2$.

\noindent \\
{\bf (3.1) Example:} In the simplest case $m_1=m_2=1$ the irregular polarizations are
of the form $\L=(\k/n,\, 1-\k/n,\, 1/n)$, $0 \le \k \le n$. The set of semistable morphisms
may be empty for some polarizations, for instance, if $n > \k + \text{dim}(S^{d_2}V^*)$
and $\k/n < \l_1 < (\k+1)/n$. We see this from the semistability conditions: $\f$ is in
$W^{ss}(G,\L)$ if and only if $\f$ is not equivalent to a morphism of the form
\bdm
\left[
\begin{array}{ll}
\star & 0_{\k+1,1} \\
\star & \star
\end{array}
\right] \qquad \text{or} \qquad \left[
\begin{array}{ll}
0_{n-\k,1} & \star \\
\star & \star
\end{array}
\right].
\edm
The following two conditions are sufficient to guarantee the existence of semistable
morphisms corresponding to all polarizations:
\bdm
n \le {r+d_2-1 \choose r} \quad \text{and} \quad n \le {r-1 + d_1 \choose r-1}.
\edm
To see this we choose a nonzero linear form $\psi \in V^*$ and a linear complement
$U \subset V^*$ of the subspace generated by $\psi$. We choose linearly independent
elements
\bdm
\f_{11}, \ldots, \f_{n1} \in S^{d_1} U \quad \text{and} \quad \psi_{12}, \ldots, \psi_{n2}
\in S^{d_2-1}V^*.
\edm
We put $\f_{i2} = \psi_{i2} \, \psi$ for $1 \le i \le n$. We claim that $\f$ is not equivalent
to a matrix having a zero entry. Indeed, the entries from the second column are linearly
independent, and no linear combination of entries from the second column can divide
a linear combination of entries from the first column, because the former is divisible
by $\psi$, whereas the latter is not.

Under the above conditions on $n$ we see that the singular polarizations are precisely the
irregular ones. Indeed, for $\l_1=\k/n$ there exist properly semistable morphisms: just choose
$\f_{11}, \ldots, \f_{\k 1}$ linearly independent in $S^{d_1}U$, put $\f_{i1}=0$ for
$\k +1 \le i \le n$ and $\f_{i2}= \psi_{i2} \psi$ as above. \\

We now turn to the embedding into the action of the reductive group. In order to
apply the theory from section 2 we need to assume a priori that $\a_2 > 0$, that is
$\l_2 > a \l_1$. According to (2.8), a point $(\xi,\g) \in \Wtilda$ is semistable if and only if
\bdm
(\xi,\g) \nsim (\xi_{ki},\g_{lk}) \quad \text{with} \quad \frac{l}{n} > k\l_1 + i \a_2.
\edm
Here $\x_{ki}$, $\g_{lk}$ are matrices of the form
\bdm
\x_{ki} = \left[
\begin{array}{ll}
\star & 0_{k,m_2-i} \\
\star & \star
\end{array}
\right], \qquad \g_{lk} = \left[
\begin{array}{ll}
\star & 0_{l,p_1-k} \\
\star & \star
\end{array}
\right].
\edm
Because of the special form of $\x$ we must have $k \le m_1 + a i$. We conclude that,
in order to ensure the semistability of $(\x,\g)$, it is enough to require the condition
\bdm
\g \nsim \g_{lk} \quad \text{with}
\quad \frac{l}{n} > k\l_1 + i \a_2 \quad \text{and} \quad k \le m_1 + a i.
\edm
Actually, it is enough to require
\begin{align*}
\tag{3.2}
\g \nsim \g_{lk} \quad \text{with}
\quad \frac{l}{n} > k\l_1 + i\a_2,\ m_1 + (i-1)a < k \le m_1 + i a,\ 0 \le i \le m_2.
\end{align*} 

Before we state the main result of this section we introduce some linear algebra constants very similar in definition to the constants $c_l$ and $d_i$ from
\cite{drezet-trautmann}. Given integers $1 \le i \le m_2-1$ and $1 \le j \le m_2a$
we denote by $k(i,j)$ the maximal dimension of a vector space $U \subset M_2 \tensor
H_{12}$ which is not contained in $M_2' \tensor H_{12}$ for any subspace $M_2' 
\subset M_2$ of dimension $i$ and for which there is a subspace in $M_2^* \tensor
A_{21}$ of dimension at least $j$ orthogonal to $U$ under the canonical bilinear
map
\bdm
(M_2 \tensor H_{12}) \times (M_2^* \tensor A_{21}) \lra H_{11}.
\edm
Given an integer $2 \le i \le m_2$ and a vector space $M$ of dimension $i$ we let
$k(i)$ be the maximal dimension of a subspace $U \subset M \tensor H_{12}$
which is not contained in $M' \tensor H_{12}$ for any proper subspace
$M' \subset M$, and for which there is a nonzero subspace in $M^* \tensor A_{21}$,
orthogonal under the canonical bilinear map
\bdm
(M \tensor H_{12}) \times (M^* \tensor A_{21}) \lra H_{11}.
\edm

\noindent \\
{\bf (3.3) Claim:} \emph{Let $W$ be the space of morphisms}
\bdm
\f : m_1 \O(-d_1) \oplus m_2 \O(-d_2) \lra n\O \qquad \text{\emph{on}} \quad
\P^r=\P(V).
\edm
\emph{Assume that $m_1 < a = \text{dim}(S^{d_1-d_2}V^*)$. Let $\L$ be a nonsingular
polarization satisfying the following conditions:}
\begin{align*}
\l_2 & > a\l_1, \\
i n \l_2 & \ge nm_1 \l_1 + k(i,m_2a-i a-m_1) && \text{\emph{for}} \quad 1 \le i \le m_2-1, \\
i n \l_2 & \ge nm_1 \l_1 + k(i) && \text{\emph{for}} \quad 2 \le i \le m_2, \\
n \l_1 + i n\l_2 & \ge k(i,m_2a-i a-a+1) && \text{\emph{for}} \quad 1 \le i \le m_2-1.
\end{align*}
\emph{Then $W^{ss}(G,\L)$ admits a geometric quotient modulo $G$, which is a
quasiprojective variety.}

\noindent \\
\emph{Proof:} Let $\f$ be in $W^{ss}(G,\L)$. According to (2.6), we need to show that
$(\xi,\g(\f))$ is stable. Perturbing slightly $\L$ we can ensure that $\Ltilda$ is a nonsingular
polarization and, at the same time, that $W^{ss}(G,\L)$ has not changed.
Thus $\Wtilda^{ss}(\Gtilda,\Ltilda)=\Wtilda^s(\Gtilda,\Ltilda)$ and it is enough to show that
$\g(\f)$ satisfies condition (3.2).

We will argue by contradiction. Assume that $\g(\f)$ is equivalent to some $\g_{lk}$ with
\bdm
\frac{l}{n} > k\l_1 +i\a_2, \quad m_1+(i-1)a < k \le m_1+ia, \quad 0 \le i \le m_2.
\edm
Let $\psi$, $\psi'$, $\psi''$ denote the truncated matrices consisting of the first $l$ rows
of $\f$, $\f'$, $\f''$.
In view of the comments before (2.8), we may assume that the vector subspace in
$M_1 \tensor H_{11} \, \oplus \, M_2 \tensor H_{12}$ spanned by the rows of $\psi$
is orthogonal to a subspace, denoted $\Ker(\psi)$, inside $M_1^* \, \oplus \, M_2^* \tensor
A_{21}$, of dimension at least $p_1-k$. Orthogonality here is understood under the canonical
bilinear map
\bdm
(M_1 \tensor H_{11} \, \oplus \, M_2 \tensor H_{12}) \times (M_1^* \, \oplus \, M_2^* \tensor A_{21}) \lra H_{11}.
\edm
By analogy, considering the pairings
\bdm
(M_1 \tensor H_{11}) \times M_1^* \lra H_{11} \qquad \text{and} \qquad
(M_2\tensor H_{12}) \times (M_2^* \tensor A_{21}) \lra H_{11},
\edm
we define the orthogonal subspaces $\Ker(\psi')$ and $\Ker(\psi'')$. The dimensions of the
corresponding spaces are denoted by $\text{ker}(\psi)$, $\text{ker}(\psi')$, $\text{ker}(\psi'')$.

\noindent \\
\emph{Case $i=0,\ 1 \le k \le m_1$.} As $\text{ker}(\psi'') \ge m_2a-k > (m_2-1)a$
we see that $\Ker(\psi'')$ intersects nontrivially every copy of $A_{21}$ inside
$M_2^* \tensor A_{21} \isom k^{m_2} \tensor A_{21}$.
Thus $\psi''=0$. Replacing possibly $\f$ with an equivalent morphism we may assume
that $\text{ker}(\f') \ge m_1-k$. We obtain that $\f$ is equivalent to $\f_{k0}$, which
contradicts the semistability of $\f$.

\noindent \\
\emph{Case $i=1,\ m_1 < k < a$.} Again $\psi''=0$. As $l/n > m_1\l_1 + \a_2 > m_1 \l_1$ we arrive at the contradiction $\f \sim \f_{m_1,0}$.

\noindent \\
\emph{Case $i=1,\ a \le k \le m_1+a$.} If $\psi'' \sim \psi_1''$, then $\text{ker}(\psi'')
=(m_2-1)a$ hence, after possibly replacing $\f$ with an equivalent morphism,
we may assume that
\bdm
\text{ker}(\psi') \ge m_1+m_2a-k-(m_2-1)a=m_1+a-k.
\edm
Moreover, from $l/n > k\l_1 + \a_2 = (k-a)\l_1+\l_2$ we get $l \ge l_{k-a,1}$.
Thus $\f \sim \f_{k-a,1}$, contradiction.

If $\psi'' \nsim \psi_1''$, then the rows of $\psi''$ span a space of dimension at most
\bdm
k(1,m_2a-k) \le k(1,m_2a-a-m_1).
\edm
From the hypotheses of the claim we have
\bdm
l > n k\l_1 + n\a_2 \ge na\l_1 + n\a_2 = n\l_2 \ge nm_1 \l_1 + k(1,m_2a-k)
\edm
forcing $l \ge l_{m_1,0}+ k(1,m_2a-k)$. We conclude that $\f \sim \f_{m_1,0}$, contradiction.

\noindent \\
\emph{Case $i \ge 2,\ m_1+(i-1)a < k < ia$.} If $\psi'' \sim \psi_{i-1}''$, then, taking into
account the inequalities
\bdm
\frac{l}{n}> m_1\l_1 + (i-1)a\l_1 + i\a_2 > m_1 \l_1 + (i-1)\l_2,
\edm
we obtain the contradiction $\f \sim \f_{m_1,i-1}$. Assume now that $\psi'' \nsim \psi_{i-1}''$. The rows of $\psi''$ must span a vector space of dimension at most
$k(i-1,m_2a-k)$. From the hypotheses of the claim we have the inequalities
\begin{align*}
l & > n k \l_1 + ni\a_2 \ge n(m_1+(i-1)a+1)\l_1+ni\a_2 \\
& > nm_1\l_1 + n\l_1 + n(i-1)\l_2 \\
& \ge nm_1\l_1 + k(i-1, m_2a-ia+1) \\
& \ge nm_1\l_1 + k(i-1,m_2a-k).
\end{align*}
From the above we get $l \ge l_{m_1,0} + k(i-1,m_2a-k)$, leading to the contradiction
$\f \sim \f_{m_1,0}$.

\noindent \\
\emph{Case $i \ge 2,\ ia \le k \le m_1+ia$.} As above, if $\psi'' \sim \psi_{i-1}''$
we get a contradiction. Assume now that $\psi'' \sim \psi_i''$ and $\psi'' \nsim \psi_{i-1}''$.
If $\text{ker}(\psi'') = (m_2-i)a$, then, replacing possibly $\f$ with an equivalent morphism,
we may assume that
\bdm
\text{ker}(\psi') \ge m_1+m_2a-k-(m_2-i)a = m_1+ia-k.
\edm
Moreover, from the inequalities $l/n > k\l_1 + i \a_2 = (k-ia)\l_1 + i\l_2$
we get $l \ge l_{k-ia,i}$. These lead to the contradiction $\f \sim \f_{k-ia,i}$.

If $\text{ker}(\psi'') > (m_2-i)a$, then the rows of $\psi''$ span a space of dimension
at most $k(i)$. From the hypotheses of the claim we have the inequalities
\bdm
l > nk\l_1 + ni\a_2 \ge nia\l_1 + ni\a_2 = n i\l_2 \ge nm_1\l_1 + k(i).
\edm
Thus $l \ge l_{m_1,0}+k(i)$ leading to the contradiction $\f \sim \f_{m_1,0}$.

Finally, assume that $\psi'' \nsim \psi_i''$. Then the rows of $\psi''$ span a space of
dimension at most $k(i,m_2a-k)$. From the hypotheses of the claim we have the
inequalities
\begin{align*}
l & \ge nk\l_1 + ni\a_2 \ge nia\l_1 + ni\a_2 = ni \l_2 \\
& \ge nm_1\l_1 + k(i,m_2a-ia-m_1) \\
& \ge nm_1\l_1 + k(i,m_2a-k).
\end{align*}
From the above we get $l \ge l_{m_1,0}+ k(i,m_2a-k)$, hence the contradiction
$\f \sim \f_{m_1,0}$.

\noindent \\
The disadvantage of (3.3) is that the linear algebra constants $k(i,j)$ and $k(i)$
are difficult to compute in general. We carry out the computations only for a particular
kind of morphisms in section 7. Sections 4, 5, 6 and 8 are not direct applications of (3.3):
in the above proof a lot of information is lost in the course of estimating $l$. We can get
more precise statements by examining each case separately, yet the arguments will be
just reworkings of the arguments from the proof of (3.3)


\section{Morphisms of the Form $m\O(-d_1)\oplus 2\O(-d_2) \lra n\O$}

We fix integers $d_1 > d_2 >0$, we fix a vector space $V$ of dimension $r+1$
and we consider morphisms
\bdm
\f = (\f',\f'') : m\O(-d_1) \oplus 2\O(-d_2) \lra n\O \qquad \text{on} \quad \P^r = \P(V).
\edm
The singular values for $\l_1$ are among those values for which there are integers
$0 \le \k \le n$ and $0 \le p \le m$ such that
\bdm
\frac{\k}{n} = p\l_1 \quad \text{or} \quad \frac{\k}{n}=p\l_1 + \l_2 \quad \text{or} \quad
\frac{\k}{n} = p\l_1 + 2\l_2.
\edm
If both $m$ and $n$ are even we can choose $\k=n/2$, $p=m/2$ and the second equation
will be satisfied for all $\l_1$, in other words all polarizations will be irregular.
To avoid this, we will assume in the sequel that either $m$ or $n$ is odd.
Under this assumption the singular values for $\l_1$ are among the numbers $\k/pn$,
$0 \le \k \le n$, $1 \le p \le m$.

The semistability conditions in the special case $m=1$ read as follows:
for $\k/n < \l_1 < (\k+1)/n$ the morphism $\f$ is semistable if and only if $\f$ is not
equivalent to a matrix having one of the following forms:
\bdm
\f_1 = \left[
\begin{array}{ll}
\star & 0_{l_1,2} \\
\star & \star
\end{array}
\right], \quad \f_2 = \left[
\begin{array}{ll}
0_{l_2,2} & \star \\
\star & \star
\end{array}
\right], \quad \f_3 = \left[
\begin{array}{ll}
\star & 0_{l_3,1} \\
\star & \star
\end{array}
\right], \quad \f_4 = \left[
\begin{array}{ll}
0_{l_4,1} & \star \\
\star & \star
\end{array}
\right],
\edm
where
\bdm
l_1=\k+1,\qquad l_2= \left[ \frac{n-\k+1}{2} \right], \qquad l_3 = \left[ \frac{n+\k}{2} \right] +1,
\qquad l_4=n-\k.
\edm
The set of semistable morphisms may be empty, for instance when $n > \k + 2 \, \text{dim}(S^{d_2}V^*)$. However, the following three conditions are enough to guarantee the existence
of semistable morphisms corresponding to all chambers (we are still in the case $m=1$):
\bdm
\tag{4.1}
n \le \text{dim}(S^{d_2-1}V^*) + l_3 -1, \quad n \le 2\, \text{dim}(S^{d_2-1}V^*)+l_1-1,
\quad n \le \text{dim}(S^{d_1} U).
\edm
We refer to (3.1) for the meaning of $U$. Explicitly, the above conditions take the form
\begin{align*}
n & \le 2\, {r+d_2-1 \choose r} + \k -1 && \text{when} \quad n+\k \quad \text{is odd}, \\
n & \le 2\, {r+d_2-1 \choose r} + \k     && \text{when} \quad n+\k \quad \text{is even}, \\
n & \le \ \ {r+d_1-1 \choose r-1}.
\end{align*}

We are not able to say precisely what are the singular values for $\l_1$ in general,
however, in the special case $m=1$, and under the assumption (4.1),
one can see as at (3.1) that the singular values are $\k/n$, $0 \le \k \le n$.

We now turn to the embedding into the action of the reductive group.
In order to apply the theory from section 2 we need to assume a priori that
$\a_2 > 0$, that is $\l_1 < 1/(2a+m)$. According to (3.2), a point $(\x,\g) \in \Wtilda$
is semistable if $\g \nsim \g_{lk}$
\begin{align*}
\text{with} \quad \frac{l}{n} & > k\l_1 && \text{and} \quad 0 \le k \le m, \\
\text{or with} \quad \frac{l}{n} & > k\l_1 + \a_2 && \text{and} \quad m+1 \le k \le a+m, \\
\text{or with} \quad \frac{l}{n} & > k\l_1 + 2\a_2 && \text{and} \quad a+m+1 \le k \le 2a+m-1.
\end{align*}

\noindent \\
{\bf (4.2) Claim:} \emph{Let $m$ and $n$ be positive integers at least one of which is odd.
Let $W$ be the space of morphisms}
\bdm
\f: m\O(-d_1) \oplus 2\O(-d_2) \lra n\O \qquad \text{\emph{on}} \quad \P^r=\P(V).
\edm
\emph{Let $0 < \l_1 < 1/(2a+m)$ be a nonsingular value. Assume that either the conditions}
\bdm
\tag{i}
m< {r-1+d_1-d_2 \choose r-1} \qquad \text{\emph{and}} \qquad
\l_1 \le \frac{n-\text{\emph{dim}}(S^{d_2-1}V^*)}{(a+m-1)n}
\edm
\emph{or the conditions}
\bdm
\tag{ii}
m< {r+d_1-d_2 \choose r} \qquad \text{\emph{and}} \qquad
\l_1 \le \frac{n-2 \, \text{\emph{dim}}(S^{d_2-1}V^*)}{3mn}
\edm
\emph{are satisfied. Then the set of semistable morphisms admits a geometric quotient
$W^{ss}(G,\L)/G$, which is a quasiprojective variety.}

\noindent \\
\emph{Proof:} Let $\f$ be in $W^{ss}(G,\L)$. According to (2.6), we need to show that
$(\x,\g(\f))$ is semistable. We argue by contradiction. 
Assume that $\g(\f) \sim \g_{lk}$ with $l/n > k\l_1$ and $0 \le k \le m$.
Let $\psi = (\psi',\psi'')$ denote the truncated matrix consisting of the first $l$ rows of $\f$.
By assumption $\psi$ has kernel inside $k^m \oplus S^{d_1-d_2}V^* \oplus S^{d_1-d_2}V^*$
of dimension at least $2a+m-k$ which is greater than $a+m$ because, by hypothesis,
$m< a$. This shows that the kernel of $\psi$ intersects each copy of $S^{d_1-d_2}V^*$
nontrivially, forcing $\psi''=0$. Moreover, there are at least $m-k$ linearly independent
elements in the kernel of $\psi'$ viewed as a subspace of $k^m$. We get $\f \sim \f_{k,0}$,
which contradicts the semistability of $\f$.

Assume now that $\g(\f) \sim \g_{lk}$ with $l/n > k\l_1 + \a_2$ and $m+1 \le k \le a+m$.
Note that automatically $l \ge l_{m,0}$, thus excluding those $\g_{lk}$ with $k<a$ because,
as we saw above, the condition $k < a$ forces $\psi''=0$, yielding the contradiction
$\f \sim \f_{m,0}$.

Assume that $a \le k \le a+m$. We have $\text{ker}(\psi) \ge 2a+m-k > m$, which forces
$\text{ker}(\psi'') \ge 1$. Let $(f,g)$ be a nonzero vector of $\Ker(\psi'')$ regarded as a
subspace of $S^{d_1-d_2}V^* \oplus S^{d_1-d_2}V^*$. Assume that $f,g$ are linearly
dependent. Replacing possibly $\f$ with an equivalent morphism we may assume that
$g=0$, $f \neq 0$ so the first column of $\psi''$ is zero.
The second column of $\psi''$ is not zero because $\f \nsim \f_{m,0}$.
Now we have $\text{ker}(\psi'')=a$ and, replacing possibly $\f$ with an equivalent morphism,
we may assume that $\text{ker}(\psi') \ge a+m-k$. As $l/n > (k-a)\l_1+\l_2$ we arrive at
$\f \sim \f_{k-a,1}$, contradiction.

Assume now that $f,g$ are linearly independent. We write $f=h f_1$, $g=h g_1$ with
$f_1$, $g_1$ relatively prime and $\text{max}\{ 0,d_1-2d_2 \} \le d = \text{deg}(h)< d_1-d_2$.
The rows of $\psi''$ are of the form $(-g_1 u,f_1 u)$, hence they are vectors in a space
of dimension equal to $\text{dim}(S^{2d_2+d-d_1}V^*)$. As $\psi'' \neq 0$ we see that
the kernel of $\psi''$ consists of vectors of the form $(v f_1,v g_1)$, hence
$\text{ker}(\psi'') = \text{dim}(S^dV^*)$. Writing $b= \text{dim}(S^{d_1-d_2-1}V^*)$ we have
\bdm
m \ge \text{ker}(\psi') \ge 2a+m-k- \text{dim}(S^dV^*) \ge 2a+m-k-b
\edm
hence
\bdm
b+k \ge 2a \quad \text{so} \quad b+a+m \ge 2a \quad \text{giving} \quad m \ge a-b.
\edm
This contradicts hypothesis (i) from the statement of the claim.
Under hypothesis (ii), we would get a contradiction if we could show that the inequality
$l/n > k\l_1 + \a_2$ implies the inequality
\bdm
l \ge l_{m,0} + \text{dim}(S^{2d_2+d-d_1}V^*).
\edm
Indeed, we would arrive at the contradiction $\f \sim \f_{m,0}$. Thus, for all $a \le k \le a+m$
and $d$ we need the inequality
\bdm
n k \l_1 + n \a_2 \ge mn \l_1 + \text{dim}(S^{2d_2+d-d_1}V^*).
\edm
This would follow from the inequality
\bdm
na\l_1 + n\a_2 \ge mn\l_1 + \text{dim}(S^{d_2-1}V^*).
\edm
But the above is equivalent to the condition on $\l_1$ from hypothesis (ii).

We are left to examine the situation $\g(\f) \sim \g_{lk}$ with $l/n > k\l_1 + 2\a_2$
and $a+m+1 \le k \le 2a+m-1$. As $k\l_1 + 2\a_2 \ge (m+1)\l_1 + \l_2 + \a_2$,
we see that $l \ge l_{m,1}$. If $k < 2a$ then, as above, we get $\text{ker}(\psi'') \ge 1$.
Let $f$, $g$ be as above. If $f$, $g$ are linearly dependent, we get the contradiction
$\f \sim \f_{m,1}$. If $f$, $g$ are linearly independent, we get a contradiction as above
under hypothesis (ii). Under hypothesis (i) we would also get a contradiction if we could
show that the inequality $l/n > k\l_1 + 2\a_2$ implies the inequality
\bdm
l \ge l_{m,0} + \text{dim}(S^{2d_2+d-d_1}V^*).
\edm
Thus we need the estimate
\bdm
n(a+m+1) \l_1 + 2n\a_2 \ge mn\l_1 + \text{dim}(S^{d_2-1}V^*).
\edm
Using the relations $\a_2 = \l_2 -a\l_1$ and $2\l_2= 1-m\l_1$ we see that the above is
equivalent to the estimate on $\l_1$ from hypothesis (i).

Finally, assume that $2a \le k \le 2a+m-1$. If $\text{ker}(\psi'') \ge 1$ we get a contradiction
as above. If $\text{ker}(\psi'')=0$, then $\text{ker}(\psi') \ge 2a+m-k$. As $l/n > (k-2a)\l_1 + 2\l_2$ we arrive at the contradiction $\f \sim \f_{k-2a,2}$. This finishes the proof of the claim. \\

According to corollary 7.2.2 in \cite{drezet-trautmann}, if $\a_2 > 0$ and if $\l_2 \ge a_{21} c_1(2)/n$, then the conclusion of (4.2) is true. The constant $c_1(2)$ can be computed as
in the proof of lemma 9.1.2 (loc. cit.) One has
\bdm
c_1(2) = \frac{\text{dim}(S^{d_2}V)}{\text{dim}(S^{d_1-d_2}V)}.
\edm
Thus, according to Dr\'ezet and Trautmann, the conclusion of claim (4.2) holds if
\bdm
\l_2 \ge \frac{1}{n} {r+ d_2 \choose r}, \quad \text{that is} \quad
\l_1 \le \frac{1}{m} \left( 1 - \frac{2}{n} {r+d_2 \choose r} \right).
\edm
Our result is not contained in Dr\'ezet and Trautmann's result. \\

In the special case $m=1$, $r \ge 2$ the estimates on $m$ from (4.2)(i) and (ii) are
automatically fulfilled, so we have:

\noindent \\
{\bf (4.3) Claim:} \emph{Let $W$ be the space of morphisms}
\bdm
\f: \O(-d_1) \oplus 2\O(-d_2) \lra n\O \qquad \text{\emph{on}} \quad \P^r, \quad r \ge 2.
\edm
\emph{Let $0 < \l_1 < 1/(2a+1)$ be a nonsingular value satisfying one of the following two
conditions:}
\bdm
\l_1 \le \frac{1}{a}- \frac{1}{an} {r+d_2-1 \choose r} \qquad \text{\emph{or}} \qquad
\l_1 \le \frac{1}{3} - \frac{2}{3n} {r+d_2-1 \choose r}.
\edm
\emph{Then the set of semistable morphisms admits a geometric quotient modulo $G$,
which is a quasiprojective variety.} \\

At the end of this section we would like to spell out a simple case in which we can
say for sure that the above claim is nonvacuous, that is, in which we know that
$W^{ss}(G,\L)$ is not empty. The case of the left-most chamber $0 < \l_1 < 1/n$
is completely understood, cf. lemma 9.3.1 in \cite{drezet-trautmann}. Let us take
$1/n < \l_1 < 2/n$ and
\bdm
\tag{4.4}
n=2 {r+d_2-1 \choose r} + 1, \qquad d_1 \le 2d_2-2.
\edm
Condition (4.1) is satisfied for $\k=1$ because $n+\k$ is even. The conditions on
$\l_1$ from (4.3) reduce to $0 < \l_1 < 1/(2a+1)$. We have $1/n < 1/(2a+1)$ because
of the second inequality in (4.4). We conclude that claim (4.3) is nonvacuous for morphisms
satisfying (4.4).


\section{Morphisms of the Form $\O(-d-1) \oplus 3\O(-d) \lra n\O$}

We fix an integer $d>0$, we fix a vector space $V$ over $k$ of dimension 3 and we consider morphisms
\bdm
\f=(\f',\f''): \O(-d-1) \oplus 3\O(-d) \lra n\O \qquad \text{on} \quad \P^2=\P(V).
\edm
Keeping the notations from section 2 we have:
\bdm
a=3,\ m_1=1,\ m_2=3,\ p_1=10,\ \a_2= \l_2-3\l_1 = \frac{1-10\l_1}{3}.
\edm
The singular values for $\l_1$ are among those values for which there are integers
$0 \le \k \le n$ and $0 \le p \le 3$ such that
\bdm
\frac{\k}{n}=p \l_2 \qquad \text{or} \qquad \frac{\k}{n}= \l_1+p \l_2.
\edm
Using the relation $\l_1 + 3\l_2 =1$ we see that the singular values for $\l_1$
are among the numbers $\k/2n$, $0 \le \k \le 2n$.

We write $\f$ as a matrix $(\f_{ij})_{1\le i \le n,\, 1 \le j \le 4}$ with $\f_{i1} \in S^{d+1}V^*$
while $\f_{ij} \in S^dV^*$ for $j=2,\, 3,\, 4$. The morphism $\f$ is semistable if and only if
it is not equivalent to a matrix having one of the following forms:
\bdm
\f_1= \left[
\begin{array}{ll}
\star & 0_{l_1,3} \\
\star & \star
\end{array}
\right], \qquad \f_2 = \left[
\begin{array}{ll}
0_{l_2,3} & \star \\
\star & \star
\end{array}
\right], \qquad \f_3 = \left[
\begin{array}{ll}
\star & 0_{l_3,2} \\
\star & \star
\end{array}
\right],
\edm
\bdm
\f_4= \left[
\begin{array}{ll}
0_{l_4,2} & \star \\
\star & \star
\end{array}
\right], \qquad \f_5 = \left[
\begin{array}{ll}
\star & 0_{l_5,1} \\
\star & \star
\end{array}
\right], \qquad \f_6 = \left[
\begin{array}{ll}
0_{l_6,1} & \star \\
\star & \star
\end{array}
\right].
\edm
Here each $\f_i$ has a zero submatrix with $l_i$ rows. The integers $l$ are
the smallest integers satisfying
\bdm
\frac{l_1}{n} > \l_1, \quad
\frac{l_2}{n} > \l_2, \quad
\frac{l_3}{n} > \l_1+\l_2, \quad
\frac{l_4}{n} > 2\l_2, \quad
\frac{l_5}{n} > \l_1+2\l_2, \quad
\frac{l_6}{n} > 3\l_2.
\edm

In order to apply the theory from section 2 we need to have $\a_2 > 0$, that is we need to assume a priori that $\l_1 < 1/10$. According to (3.2), in order to show that $(\x,\g)$
is semistable, it is enough to show that $\g \nsim \g_{lk}$
\begin{align*}
\text{with} \quad \frac{l}{n} & > \l_1 && \text{and} \quad k=1, \\
\text{or with} \quad \frac{l}{n} & > k\l_1 + \a_2 && \text{and} \quad k=2,\, 3,\, 4, \\
\text{or with} \quad \frac{l}{n} & > k\l_1 + 2\a_2 && \text{and} \quad k=5,\, 6,\, 7, \\
\text{or with} \quad \frac{l}{n} & > k\l_1 + 3\a_2 && \text{and} \quad k=8,\, 9.
\end{align*}

\noindent \\
{\bf (5.1) Claim:} \emph{Let $W$ be the space of morphisms of sheaves on $\P^2$
of the form}
\bdm
\f : \O(-d-1) \oplus 3\O(-d) \lra n\O.
\edm
\emph{Then for any nonsingular value for $\l_1$ satisfying the conditions}
\begin{align*}
0 < \l_1 & < \frac{1}{10}, \\
\l_1 & \le \frac{2}{5} - \frac{3}{10n} (d^2+d), \\
\l_1 & \le 1 - \frac{3}{4n}(d^2+3d), \\
\l_1 & \ge -\frac{1}{2} + \frac{3}{4n}(d^2+d), \\
\l_1 & \ge -2 + \frac{3}{n}(d^2+2d)
\end{align*}
\emph{the set of semistable morphisms admits a geometric quotient modulo $G$,
which is a quasiprojective variety.}

\emph{Notice that for $d$ and $n$ large, and for $n$ of order $3d^2/2$ the last four
inequalities follow from the first.}

\noindent \\
\emph{Proof:} Let $\f$ be in $W^{ss}(G,\L)$. According to (2.6), we need to show that
$(\x,\g(\f))$ is semistable. We argue by contradiction. Assume that $\g(\f) \sim \g_{lk}$
with $l/n > \l_1$ and $k=1,\, 2$. Note that $l \ge l_1$. Let $\psi = (\psi',\psi'')$
denote the truncated matrix consisting of the first $l$ rows of $\f$.
By assumption $\psi$ has kernel inside $k \oplus V^* \oplus V^* \oplus V^*$ of
dimension at least $10-k$. For dimension reasons $\Ker(\psi)$ intersects each copy of
$V^*$ nontrivially, forcing $\psi''=0$. 
We get $\f \sim \f_1$, which contradicts the semistability of $\f$.
The same argument also works in the case
$k=3$ and $\text{ker}(\psi'')=7$.

Assume now that $k=3$, $l/n > 3\l_1 + \a_2 = \l_2$ and $\text{ker}(\psi'')=6$.
Replacing $\f$ with an equivalent morphism we may assume that $\psi'=0$.
From remark (5.2) below we see that two columns of $\psi''$ vanish.
As $l \ge l_2$ we get $\f \sim \f_2$, contradiction.

Assume that $k=4$ and $l/n > 4\l_1 + \a_2 = \l_1 + \l_2$. As $\text{ker}(\psi'') \ge 5$
and $l \ge l_3$, we can apply (5.2) and we get the contradiction $\f \sim \f_3$.
The same argument works if $k=5$ and $\text{ker}(\psi'')\ge 5$.

Assume now that $k=5$, $l/n > 5\l_1 + 2\a_2=(2-5\l_1)/3$ and that $\text{ker}(\psi'')=4$.
The matrix $\psi''$ has the form given at remark (5.3) below, so its rows are elements
in a vector space of dimension at most equal to the dimension of $S^{d-1}V^*$.
If we could show that $l > \text{dim}(S^{d-1}V^*)$, then we would conclude that $\f$
has a zero row, which would be a contradiction. Thus we need the inequality
\bdm
\frac{n(2-5\l_1)}{3} \, \ge \, {d+1 \choose 2}.
\edm
But this is equivalent to the second condition on $\l_1$ from the statement of the claim.

Assume that $k=6$ and $l/n > 6\l_1 +2\a_2 = 2\l_2$. Thus $l \ge l_4$ and, from the above inequality, $l \ge l_1 + \text{dim}(S^{d-1}V^*)$. If $\text{ker}(\psi'')=4$, then $\psi''$
has the form given at remark (5.3) below and we arrive at the contradiction
$\f \sim \f_1$. If $\text{ker}(\psi'')=3$, then we may assume that $\psi'=0$.
Let $\eta$ be a $3\times 3$-matrix with entries in $V^*$ whose columns are
linearly independent vectors in $\Ker(\psi'')$. Each column of $\eta$ must contain
at least two linearly independent elements, otherwise we get the contradiction
$\f \sim \f_4$. Also, each row of $\eta$ must contain at least two linearly independent
elements. Indeed, if
\bdm
\eta = \left[
\begin{array}{ccc}
u_1 & u_2 & u_3 \\
v_1 & v_2 & v_3 \\
0 & 0 & w_3
\end{array}
\right]
\edm
with $u_1$, $v_1$ linearly independent, $u_2$, $v_2$ linearly independent, then
\bdm
\left|
\begin{array}{cc}
u_1 & u_2 \\
v_1 & v_2
\end{array}
\right| \neq 0
\edm
and we conclude that the first two columns of $\psi''$ are zero, again a contradiction.

Thus $\eta$ satisfies the hypotheses of remark (5.4) from below. We must have
$\eta=\eta_1$ or $\eta=\eta_2$. In the first case the rows of $\psi''$ are elements in
a vector space of dimension at most equal to the dimension of $S^{d-1}V^*$.
As $l > \text{dim}(S^{d-1}V^*)$, we conclude that $\f$ has a zero row, contradiction.

Assume now that $\eta = \eta_2$. From remark (5.5) we know that the rows of
$\psi''$ are elements in a vector space of dimension at most $(d^2+3d)/2$.
The third condition on $\l_1$ from the statement of the claim is equivalent to the
inequality $2n\l_2 \ge (d^2+3d)/2$. This shows that $l > (d^2+3d)/2$, forcing $\f$
to have a zero row, again a contradiction.

The case $k=7$ follows from the last two conditions on $\l_1$ from the claim and the cases $k=8$ and $k=9$ are analogous. This finishes the proof of the claim.

\noindent \\
{\bf (5.2) Remark:} Let $\psi$ be an $l \times 3$-matrix with entries in $S^dV^*$.
We consider the space of $3 \times1$-matrices $\eta$ with entries in $V^*$
such that $\psi \eta = 0$ and we call it $\Ker(\psi)$. Let $\text{ker}(\psi)$ denote
its dimension. If $\text{ker}(\psi) \ge 5$, then two columns of $\psi$ must vanish.

\noindent \\
{\bf (5.3) Remark:} With the notations from (5.2), if $\text{ker}(\psi)=4$,
then $\psi$ has the form
\bdm
\left[
\begin{array}{ccc}
0 & uf_1 & vf_1 \\
\vdots & \vdots & \vdots \\
0 & uf_l & vf_l
\end{array}
\right]
\edm
with $u,\, v \in V^*$ linearly independent one-forms and $f_i \in S^{d-1} V^*$ for
$1 \le i \le l$.

\noindent \\
Next we quote 5.6 from \cite{maican}:

\noindent \\
{\bf (5.4) Remark:} Let $\eta$ be a $3 \times 3$-matrix with entries in $V^*$.
Assume that $\text{det}(\psi)=0$ and that $\eta$ is not equivalent, modulo elementary operations on rows or on columns, to a matrix
having a zero row, or having a zero column, or having a zero $2 \times 2$-submatrix.
Then $\eta$ is equivalent to one of the following matrices:
\bdm
\eta_1 = \left[
\begin{array}{crc}
X & Y & 0 \\
Z & 0 & Y \\
0 & -Z & X
\end{array}
\right] \quad \text{or} \quad \eta_2 = \left[
\begin{array}{lll}
X & Y & Z \\
Y & a_1 X + a_2 Y & a_3 X + a_4 Y + a_5 Z \\
Z & a_6 X + a_7 Y + a_8 Z & a_9 X + a_{10} Z
\end{array}
\right]
\edm
with nonzero constants $a_1, \ldots, a_{10}$ in the ground field $k$.
Here $\{ X,\, Y,\, Z \}$ is a basis of $V^*$.

\noindent \\
{\bf (5.5) Remark:} The kernel of $\eta_2$ inside $S^d V^* \oplus S^d V^* \oplus S^d V^*$
is of dimension at most $(d^2+3d)/2$.

\noindent \\
\emph{Proof:} The elements of $\Ker(\eta_2)$ are of the form
\bdm
f(-Y,X,0) + g(-Z,0,X)+h(0,-Z,Y)
\edm
with $(f,g,h)$ determined modulo multiples of $(Z,-Y,X)$.
Without loss of generality we may assume that $h$ depends only on $Y$ and $Z$.
We have
\begin{align*}
0 = f(-Y^2+a_1X^2+a_2XY)+g(-Y\!Z+a_6 X^2+ a_7 XY + a_8 X\!Z) \\
 + h(-a_1X\!Z-a_2Y\!Z+a_6 XY+ a_7 Y^2 + a_8 Y\!Z). \qquad \qquad
\end{align*}
This shows that $f$ is uniquely determined by $g$ and $h$.
Hence $\Ker(\eta_2)$ is of dimension at most
\bdm
{d+1 \choose 2} + d = \frac{d^2+3d}{2}.
\edm


\section{Morphisms of the Form $m\O(-d-1) \oplus 3\O(-1) \lra n\O$}

We fix an integer $d > 0$, we fix a vector space $V$ of dimension 3, and we consider morphisms
\bdm
\f = (\f',\f'') : m\O(-d-1) \oplus 3\O(-1) \lra n\O \qquad \text{on} \quad \P^2 = \P(V).
\edm
We write $b= \text{dim}(S^{d-1}V^*)$ and, using the notations from section 2, we have
\bdm
m_1 = m,\ \ m_2=3,\ \ p_1=3a+m,\ \ \a_2 = \l_2 - a\l_1 = \frac{1-(3a+m)\l_1}{3}.
\edm
The singular values for $\l_1$ are among those values for which there are integers
$0 \le \k \le n$, $0 \le p \le m$ and $0 \le q \le 3$ such that $\k/n = p\l_1 + q\l_2$.
If both $m$ and $n$ are divisible by 3 we can take $\k = n/3$, $p=m/3$, $q=1$
and we see that all values for $\l_1$ are irregular. If either $m$ or $n$ is not divisible by 3, which will be our assumption in the sequel, then, using the relation $m\l_1 + 3\l_2=1$,
we see that the singular values for $\l_1$ are among the numbers $\k/pn$ with
$0 \le \k \le 2n$, $1 \le p \le 2m$.

In order to apply the theory from section 2 we need to assume a priori that $a_2>0$,
that is $\l_1 < 1/(3a+m)$. According to (3.2), in order to show that $(\x,\g)$ is semistable,
it is enough to show that $\g \nsim \g_{lk}$
\begin{align*}
\text{with} \quad \frac{l}{n} & > k\l_1 && \text{and} \quad 0 \le k \le m, \\
\text{or with} \quad \frac{l}{n} & > k\l_1 + \a_2 && \text{and} \quad m+1 \le k \le a+m, \\
\text{or with} \quad \frac{l}{n} & > k\l_1+2\a_2 && \text{and} \quad a+m+1 \le k \le 2a+m, \\
\text{or with} \quad \frac{l}{n} & > k\l_1 + 3\a_2 && \text{and} \quad 2a+m+1 \le k \le 3a+m-1.
\end{align*}

\noindent \\
{\bf (6.1) Claim:} \emph{Let $m$ and $n$ be positive integers one of which is not divisible by 3. Let $W$ be the space of morphisms}
\bdm
\f : m\O(-d-1) \oplus 3\O(-1) \lra n\O \qquad \text{\emph{on}} \quad \P^2=\P(V).
\edm
\emph{Assume that $m<a$. Then for any nonsingular value $\l_1$ satisfying}
\bdm
0 < \l_1 < \frac{1}{3a+m}, \qquad \l_1(4m-3a+3b) \le \frac{n-3}{n}, \qquad \l_1 \le
\frac{n-6}{mn}
\edm
\emph{the set of semistable morphisms admits a geometric quotient $W^{ss}(G,\L)/G$ which is a quasiprojective variety.}

\noindent \\
\emph{Proof:} Let $\f$ be in $W^{ss}(G,\L)$. The first part of the proof proceeds as at (5.1).
So let us assume that $\g(\f) \sim \g_{lk}$ with $l/n > k\l_1 + \a_2$ and $a \le k \le a+m$.
We have $\text{ker}(\psi'') \ge a+1$. Let $r$ be the dimension of the vector space spanned
by the rows of $\psi''$. According to (5.2) we would have $\text{ker}(\psi'') \le a$ if $r \ge 5$,
which is not the case. According to (5.3) we would have $\text{ker}(\psi'')\le \text{dim}(S^{d-1}V^*)$ if $r=4$, which is not the case either. Thus $r \le 3$.

Assume that $r=3$. Let $\eta$ be a $3 \times 3$-matrix with entries in $V^*$
formed from three linearly independent rows of $\psi''$.
We cannot have
\bdm
\eta = \left[
\begin{array}{ccc}
0 & 0 & u \\
v_1 & v_2 & \star \\
w_1 & w_2 & \star
\end{array}
\right] \qquad \text{with} \qquad u \neq 0, \quad \left[
\begin{array}{cc}
v_1 & v_2 \\
w_1 & w_2
\end{array}
\right] \neq 0
\edm
because this would lead to $\text{ker}(\psi'') \le \text{ker}(\eta) \le a$.
Since $\text{det}(\eta) = 0$, we can apply (5.4) to deduce that $\eta = \eta_1$, $\eta=\eta_2$ or that $\eta$ is equivalent to a matrix having a zero column. The last case leads to the
contradiction $\f \sim \f_{k-a,1}$. If $\eta = \eta_1$, we get $\text{ker}(\eta) \le \text{dim}(S^{d-1}V^*)$. If $\eta = \eta_2$, we get, in view of (5.5), the inequality $\text{ker}(\eta) \le
(d^2+3d)/2=a-1$. Both of them are contrary to $\text{ker}(\psi'') \ge a+1$.

Assume now that $r=2$. Let $\eta$ be a $2 \times 3$-matrix formed from 
two linearly independent rows of $\psi''$. If $\eta$ had two zero columns we
would arrive as before at the contradiction $\f \sim \f_{k-a,1}$.
Moreover, $\eta$ has two linearly independent elements in each row,
otherwise we would get $\text{ker}(\psi'')\le a$. If
\bdm
\eta= \left[
\begin{array}{ccc}
\star & \star & u \\
v & w & 0
\end{array}
\right]
\edm
with $u \neq 0$, then $\Ker(\psi'')$ consists of matrices of the form
\bdm
\left[
\begin{array}{rcr}
-w f_1 & \cdots & -w f_s \\
v f_1 & \cdots & v f_s \\
g_1 & \cdots & g_s
\end{array}
\right]
\edm
with each $g_j$ uniquely determined by $f_j$. This shows that $\text{ker}(\psi'')
\le \text{dim}(S^{d-1}V^*)$, contradiction. We conclude that $\eta$ is not
equivalent to a matrix having a zero entry, i.e. $\eta$ satisfies the hypothesis
of remark (6.2) from below. According to remark (6.3) we again arrive at a
contradiction: $\text{ker}(\psi'')\le a-1$.

The case $r=1$ and the remaining possibilities for $\g(\f)$ are dealt with
in a similar manner. This finishes the proof of the claim.

\noindent \\
{\bf (6.2) Remark:} Let $\eta$ be a $2 \times 3$-matrix with entries in $V^*$.
Assume that, under the canonical action of $\GL(2) \times \GL(3)$, $\eta$ is not equivalent to a matrix having a zero entry.
Then $\eta$ is equivalent to one of the following two matrices:
\bdm
\eta_3 = \left[
\begin{array}{lll}
X & Y & Z \\
Y & b_1X+b_2Y+b_3Z & c_2Y+c_3Z
\end{array}
\right]
\edm
with $b_1,\, c_2,\, c_3 \neq 0$, or
\bdm
\eta_4 = \left[
\begin{array}{lll}
X & Y & Z \\
Y & b_2Y+ b_3Z & c_1X+c_2Y+c_3Z
\end{array}
\right]
\edm
with $b_3,\, c_1 \neq 0$. Here $\{ X,Y,Z \}$ is a basis of $V^*$.

\noindent \\
Repeating the arguments from the proof of (5.5) we can estimate the dimensions
of the kernels for the above matrices:

\noindent \\
{\bf (6.3) Remark:} The kernels of $\eta_3$ and $\eta_4$ inside
$S^d V^* \oplus S^d V^* \oplus S^d V^*$ are of dimension at most
$(d^2+3d)/2$. For $d=2$ the kernel is of dimension at most 4.\\

\noindent \\
\emph{Proof:} We need to do only the case $d=2$. Keeping the notations
from (5.5) we have, say in the case $\eta=\eta_3$, the relation
\bdm
f = - \frac{(-Y\!Z+c_2XY+c_3X\!Z)g+(-b_1X\!Z-b_2Y\!Z-b_3Z^2+c_2Y^2+c_3Y\!Z)h}{-Y^2
+b_1X^2+b_2XY+b_3X\!Z}.
\edm
The kernel of $\eta_3$ can thus be parametrized by five parameters in the ground field $k$,
namely the coefficients of $g$ and $h$, because we can assume that the coefficient of
$X$ in $h$ is zero. The requirement that the above ratio be a polynomial gives algebraic
conditions in the space of parameters. For $g=Y$, $h=Y$ the numerator is not
divisible by the denominator because no monomial of the numerator is divisible by $X^2$.
Thus the algebraic conditions are nontrivial, i.e. the kernel of $\eta_3$ is
parametrized by a proper subvariety of $k^5$ and, as such, it has dimension at most 4.

The argument can be repeated for the kernel of $\eta_4$: choose $f=Y$, $h=Y$ and there
will be no possibility for $g$.

In the simplest nontrivial case $1/2mn < \l_1 < 1/(2m-1)n$ the conditions from the
claim take the form
\bdm
\frac{3a}{2m} + \frac{1}{2} < n, \qquad 5 \le n + \frac{3(a-b)}{2m}, \qquad
7 \le n.
\edm
The following conditions are also necessary to ensure the nonemptyness
of the set of semistable points:
\bdm
n < l_{m,0} + 9, \qquad n < l_{m,1} + 6, \qquad n < l_{m,2}+3.
\edm
Here $l_{m,0}=1$ while $l_{m,1}$, $l_{m,2}$ are the smallest integers satisfying
\bdm
l_{m,1} > \frac{n+1}{3}, \qquad l_{m,2} > \frac{4n+1}{6}.
\edm
These three conditions on $n$ are satisfied for $n \le 9$.
For instance, if $m=[a/2]$, all the above conditions on $n$ reduce to
$n \in \{ 7,\, 8,\, 9 \}$. It is clear that for $d$ sufficiently large the set of
semistable morphisms is not empty.


\section{Morphisms of the Form $\O(-d-2) \oplus 3\O(-d) \lra n\O$}

We fix an integer $d>0$, we fix a vector space $V$ over $k$ of dimension 3 and we
consider morphisms
\bdm
\f= (\f',\f''): \O(-d-2) \oplus 3\O(-d) \lra n\O \qquad \text{on} \quad \P^2=\P(V).
\edm
Here $m_1=1$, $m_2=3$, $a=6$ and the singular values for $\l_1$ are as in section 5.
The morphism $\f$ is semistable if and only if it is not equivalent to one of the
morphisms of the form $\f_1, \ldots, \f_6$ from section 5. From this we see that
$W^{ss}(G,\L)$ is nonempty if and only if the following conditions are satisfied:
\begin{align*}
n & < l_1 + 3 {d+2 \choose 2}, \\
n & < l_3 + 2 {d+2 \choose 2}, \\
n & < l_5 + {d+2 \choose 2}, \\
n & < l_2 + {d+4 \choose 2} + 2 {d+2 \choose 2}, \\
n & < l_4 + {d+4 \choose 2} + {d+2 \choose 2}, \\
n & < l_6 + {d+4 \choose 2}.
\end{align*}
Taking into account the definitions of $l_1, \ldots, l_6$ from section 5, the above
conditions can be rewritten as
\begin{align*}
n & \le n\l_1 + 3 {d+2 \choose 2}, \\
n & \le n\l_1 + n\l_2 + 2 {d+2 \choose 2}, \\
n & \le n\l_1 + 2n\l_2 + {d+2 \choose 2}, \\
\tag{7.1}
n & \le n\l_2 + {d+4 \choose 2} + 2{d+2 \choose 2}, \\
n & \le 2n\l_2 + {d+4 \choose 2} + {d+2 \choose 2}, \\
n & \le 3n\l_2 + {d+4 \choose 2}.
\end{align*}
Taking into account the relation $\l_1 + 3\l_2 = 1$, the first three conditions are
equivalent to
\bdm
\tag{7.2}
\l_2 \le \frac{1}{n} {d+2 \choose 2}.
\edm
The conditions on $\L$ from (3.3) for the existence of the quotient modulo $G$
read as follows:
\begin{align*}
\l_2 & > 6 \l_1, \\
n\l_2 & \ge n \l_1+ k(1,11), \\
2n\l_2 & \ge n\l_1 + k(2,5), \\
\tag{7.3}
2n\l_2 & \ge n\l_1 + k(2), \\
3n\l_2 & \ge n\l_1 + k(3), \\
n\l_1 + n\l_2 & \ge k(1,7), \\
n\l_1 + 2n\l_2 & \ge k(2,1).
\end{align*}
Next we will compute the linear algebra constants $k(i,j)$ and $k(i)$ using elementary operations
with matrices of homogeneous polynomials. We do only the most
laborious case:

\noindent \\
{\bf (7.4) Claim:} $k(2,5)= \text{\emph{dim}}(S^{d-1}V^*)= {\displaystyle {d+1 \choose 2}}$.

\noindent \\
\emph{Proof:} In the sequel $\a$ will be a matrix with 3 columns and entries in $S^d V^*$
having linearly independent rows and linearly independent columns. Also, $\b$ will be
a $3 \times 5$-matrix with entries in $S^2 V^*$ having linearly independent columns,
such that $\a \b=0$.
The constant $k(2,5)$ is the maximal number of rows that $\a$ could have for all choices
of $\a$ and $\b$.

Let us fix a basis $\{ X,Y,Z \}$ of $V^*$ and a basis $\{ u_1,\ldots, u_q \}$ of $S^{d-1}V^*$.
We consider the matrices
\bdm
\a_0 = \left[
\begin{array}{c}
u_1 \\
\vdots \\
u_q
\end{array}
\right] \left[
\begin{array}{ccc}
X & Y & Z
\end{array}
\right] \ \, \text{and} \ \, \b_0 = \left[
\begin{array}{rrrrr}
-XY & -XZ & 0 \ \, & -Y^2 & -YZ \\
X^2 & 0 \ \, & -XZ & XY & 0 \ \, \\
0 \ \, & X^2 & XY & 0 \ \, & XY
\end{array}
\right].
\edm
The choice $\a=\a_0$, $\b=\b_0$ shows that $k(2,5) \ge q$.
We will prove the converse inequality using induction on $d$,
but before that we need to make a few observations about $\a$ and $\b$.
Firstly, $\b$ cannot be equivalent to a matrix having two zeros on a column, otherwise the columns of $\a$ would be linearly dependent.
Secondly, $\b$ must have linearly independent rows. If, say, the third row of $\b$
is zero, then the matrix $\b_1$ made of the first two rows of $\b$ has
all maximal minors equal to 0. As $\b_1$ has linearly independent columns,
we see that $\b_1$ is forced to have linearly dependent rows, contradicting
observation one. Thirdly, $\a$ cannot be equivalent to a matrix having a zero entry.
If, say $\a_{13}=0$, then all maximal minors of $\b_1$ are zero, forcing $\b_1$ to have
linearly dependent columns or linearly dependent rows. This contradicts observation one,
respectively observation two.

To begin the induction assume that $d=1$. We have to show that $\a$ cannot have more
than one row. Assume that $\a$ has two rows. The third observation from above says
that $\a$ satisfies the hypotheses of remark (6.2). Then (6.3) denies the existence of $\b$,
contradiction.

Assume now that $d \ge 2$ and that $\a$ has $q+1=\text{dim}(S^{d-1}V^*)+1$ rows.
Let us put $\a'=\a \ \text{mod} \ Z$, $\b'=\b \ \text{mod} \ Z$ and let us assume that the
columns of $\b'$ are linearly independent. If $\b'$ has a zero row, then we may write
\bdm
\b= \left[
\begin{array}{cc}
\b_{11} & 0 \\
\b_{21} & \b_{22}
\end{array}
\right]
\edm
with a $1 \times 3$-matrix $\b_{11}$ that, according to observation two from above,
is not zero. As all $3 \times 3$-minors of $\b$ are zero, we get $\text{det}(\b_{22})=0$.
Moreover, $\b_{22}$ cannot have linearly dependent rows or columns, otherwise
observation one from above would be violated. Thus we may write
\bdm
\b_{22} = \left[
\begin{array}{c}
X \\ Y
\end{array}
\right] \left[
\begin{array}{cc}
u & v
\end{array}
\right]
\edm
with linearly independent $u,\, v \in V^*$. It follows that
\bdm
\a \sim \left[
\begin{array}{cc}
\a_{11} & \a_{12} \\
\a_{21} & 0
\end{array}
\right] \qquad \text{with} \qquad \a_{12} = \left[
\begin{array}{c}
v_1 \\
\vdots \\
v_q
\end{array}
\right] \left[
\begin{array}{cc}
-Y & X
\end{array}
\right]
\edm
for some $v_1, \ldots, v_q \in S^{d-1}V^*$.
This contradicts observation three from above.
Thus far we have reached the conclusion that $\b'$ cannot have a zero row.
As above, we may assume that
\bdm
\b'= \left[
\begin{array}{cc}
\b_{11}' & 0 \\
\b_{12}' & \b_{22}'
\end{array}
\right] \qquad \text{with} \qquad \b_{22}' = \left[
\begin{array}{cc}
X^2 & XY \\
XY & Y^2
\end{array}
\right]
\edm
and we get
\bdm
\a' \sim \left[
\begin{array}{cc}
\a_{11}' & \a_{12}' \\
0 & 0
\end{array}
\right] \qquad \text{with} \qquad \a_{12}' = \left[
\begin{array}{c}
v_1 \\
\vdots \\
v_d
\end{array}
\right] \left[
\begin{array}{cc}
-Y & X
\end{array}
\right]
\edm
for some homogeneous polynomials $v_1, \ldots, v_d$ of degree $d-1$ in $X$ and $Y$.
It follows that we can write
\bdm
\a = \left[
\begin{array}{r}
\a_1 \\
Z \a_2
\end{array}
\right] \quad
\text{with $\a_2$ having $q+1-d=$} {d+1 \choose 2}+1-d= {d \choose 2} +1 \text{ rows.}
\edm
From the third observation at the beginning of this proof we see that $\a_2$ satisfies the
induction hypothesis. However, $\a_2$ violates the conclusion of the induction hypothesis,
so we have arrived at a contradiction.

It remains to examine the situation in which the columns of $\b'$ are linearly dependent.
Let us write
\bdm
\b = \left[
\begin{array}{cc}
Z \b_1 & \b_2
\end{array}
\right] \qquad \text{with} \qquad \b_1 = \left[
\begin{array}{c}
u \\
v \\
w
\end{array}
\right].
\edm
By the first observation at the beginning of this proof we see that at least two among
$u$, $v$, $w$ must be linearly independent, say $u=X$, $v=Y$ and either $w=Z$ or
$w=0$. From the relation
\bdm
\a' \left[
\begin{array}{c}
X \\
Y \\
0
\end{array}
\right] = 0
\edm
we deduce, as before, the existence of $\a_2$, leading to a contradiction.
This finishes the proof of the claim.

\noindent \\
Using similar considerations, which we omit, we can compute the remaining
linear algebra constants:

\bdm
k(1,11)=0, \quad \
k(2)=k(1,7)={d+1 \choose 2}, \quad \
k(3)=k(2,1)={d+2 \choose 2} + {d+1 \choose 2}.
\edm

\noindent \\
Conditions (7.3) now take the form
\begin{align*}
\l_2 & > 6 \l_1, \\
2n \l_2 & \ge n\l_1 + {d+1 \choose 2}, \\
3n \l_2 & \ge n\l_1 + {d+2 \choose 2} + {d+1 \choose 2}, \\
n \l_1 + n\l_2 & \ge {d+1 \choose 2}, \\
n \l_1 + 2n \l_2 & \ge {d+2 \choose 2} + {d+1 \choose 2}.
\end{align*}

\noindent \\
The second inequality follows from the first and the fourth.
The third inequality follows from the first and the fifth. Using the relation
$\l_1 = 1-3\l_2$ we see that the first inequality is equivalent to $\l_2 > 6/19$.
From this and (7.2) we obtain that
\bdm
n < \frac{19}{6} {d+2 \choose 2}.
\edm
It becomes now clear that the last three inequalities in (7.1) are superfluous,
namely they are satisfied if we set $\l_2=6/19$ and $n=19(d^2+3d+2)/12$.
Eliminating all unnecessary conditions from (7.1) and (7.3), we finally arrive at:

\noindent \\
{\bf (7.5) Claim:} \emph{Let $W$ be the space of morphisms of sheaves on $\P^2$
of the form}
\bdm
\f : \O(-d-2) \oplus 3\O(-d) \lra n\O.
\edm
\emph{Then for any nonsingular polarization $\L$ satisfying}
\begin{align*}
\frac{6}{19} < \l_2 & < \frac{1}{3}, \\
\l_2 & \le \frac{1}{2} - \frac{1}{2n} {d+1 \choose 2}, \\
\l_2 & \le 1 - \frac{1}{n} {d+2 \choose 2} - \frac{1}{n} {d+1 \choose 2}, \\
\l_2 & \le \frac{1}{n} {d+2 \choose 2}
\end{align*}
\emph{the set of semistable morphisms is nonempty and admits a geometric
quotient $W^{ss}(G,\L)/G$, which is a quasiprojective variety.}

\noindent \\
\emph{Note that necessarily $n$ must satisfy the conditions}
\bdm
\frac{19}{13} {d+1 \choose 2} + \frac{19}{13} {d+2 \choose 2} < n <
\frac{19}{6} {d+2 \choose 2}.
\edm


\section{Morphisms of the Form $m\O(-d_1) \oplus \O(-d_2) \oplus \O(-d_3) \lra n\O$}

We fix integers $d_1 > d_2 > d_3 > 0$, we fix a vector space $V$ over $k$ of dimension
$r+1$ and we consider morphisms of sheaves on $\P^r=\P(V)$ of the form
\bdm
\f=(\f^1,\f^2,\f^3) : m\O(-d_1) \oplus \O(-d_2) \oplus \O(-d_3) \lra n\O.
\edm
Employing the notations from \cite{drezet-trautmann} we have:
\bdm
m_1=m, \quad m_2=1, \quad m_3 = 1, \quad p_1=m+a_{21}+a_{31}, \quad
p_2 = 1+a_{32}, \quad p_3 =1,
\edm
\bdm
m\l_1 + \l_2 + \l_3 =1,\ \a_1=\l_1, \ \a_2=\l_2-a_{21}\l_1,\ \a_3=\l_3-a_{31}\l_1
-a_{32}\l_2 +a_{32}a_{21}\l_1.
\edm
The polarization $\L=(\l_1,\l_2,\l_3,\m_1)$ is uniquely determined by the pair $(\l_1,\l_2)$
in $[0,1] \times [0,1]$. The singular polarizations are among those polarizations for which
there are integers $0 \le \k \le n$ and $0 \le p \le m$ such that
\bdm
\frac{\k}{n}=p\l_1 \qquad \text{or} \qquad \frac{\k}{n}=p\l_1+\l_2 \qquad \text{or} \qquad
\frac{\k}{n} = p\l_1+\l_3.
\edm
Equivalently, for singular polarizations the pair $(\l_1,\l_2)$ lies on one of the lines
$\l_1=\k/pn$ or $\l_2=\k/n-p\l_1$.

Given integers $0 \le \k_1 \le m$, $0 \le \k_2 \le 1$, $0 \le \k_3 \le 1$, we denote by
$l_{\k_1\k_2\k_3}$ the smallest integer satisfying
\bdm
\frac{l_{\k_1\k_2\k_3}}{n} > \k_1\l_1+\k_2\l_2+\k_3\l_3
\edm
and we consider matrices of the form
\bdm
\f^1_{\k_1\k_2\k_3} = \left[
\begin{array}{ll}
\star & 0_{l_{\k_1\k_2\k_3},m-\k_1} \\
\star & \star
\end{array}
\right], \qquad \f^i_{\k_1\k_2\k_3} = \left[
\begin{array}{l}
0_{l_{\k_1\k_2\k_3},1-\k_i} \\
\star
\end{array}
\right], \quad i=2,3.
\edm
We put $\f_{\k_1\k_2\k_3}=(\f^1_{\k_1\k_2\k_3},\f^2_{\k_1\k_2\k_3},\f^3_{\k_1\k_2\k_3})$.
According to King's Criterion of Semistability from \cite{king}, the morphism $\f$ is
semistable if and only if it is not equivalent to a morphism of the form $\f_{\k_1\k_2\k_3}$
for any choice of $\k_1$, $\k_2$, $\k_3$.

We now turn to the embedding into the action of the reductive group.
The map $\z: W \lra \Wtilda$ can be described explicitly: $\z(\f)=(\x_2,\x_3,\g(f))$ where
\bdm
\x_2 \in \text{M}_{p_1,p_2}(S^{d_1-d_2}V^*), \qquad
\x_3 \in \text{M}_{p_2,p_3}(S^{d_2-d_3}V^*), \qquad
\g(\f) \in \text{M}_{n,p_1}(S^{d_1}V^*).
\edm
Concretely, let us choose a basis $\{ X_0, \ldots, X_r \}$ of $V^*$ and let us choose bases
\begin{align*}
\{ U_i \}_{1\le i \le a_{21} }\quad & \text{of} \quad S^{d_1-d_2}V^*, \\
\{ V_j \}_{1 \le j \le a_{32}} \quad & \text{of} \quad S^{d_2-d_3}V^*, \\
\{ W_k \}_{1 \le k \le a_{31}} \quad & \text{of} \quad S^{d_1-d_3}V^*
\end{align*}
made of monomials $X_0^{i_0} \cdots X_r^{i_r}$. Then
\bdm
\x_2= \left[
\begin{array}{cc}
0 & 0 \\
U & 0 \\
0 & W
\end{array}
\right] \qquad \text{and} \qquad \x_3= \left[
\begin{array}{cccc}
0 & V_1 & \cdots & V_{a_{32}}
\end{array}
\right]^T
\edm
where
\bdm
U= \left[
\begin{array}{ccc}
U_1 & \cdots & U_{a_{21}}
\end{array}
\right]^T
\edm
while $W$ is an $a_{31} \times a_{32}$-matrix with entries $W_{kj}=W_k/V_j$ if $V_j$
divides $W_k$, otherwise $W_{kj}=0$.

In order to apply the theory from section 2 we need to assume a priori that $\a_3>0$
and $p_1 \l_1 <1$, see remark (2.7). The second condition follows from the condition
$\a_2 > 0$, which we will assume in the sequel. According to King's Criterion of
Semistability, a point $(\x_2,\x_3,\g) \in \Wtilda$ is semistable if and only if the following
conditions are satisfied:
\begin{enumerate}
\item[(i)] $(\x_2,\g) \nsim ((\x_2)_{ki}, \g_{lk})$ with $l/n > k\l_1 + i\a_2 +\a_3$,
$0 \le k \le m+a_{21}+a_{31}-1$, $1 \le i \le 1+a_{32}$. Here $\g_{lk}$ and $(\x_2)_{ki}$
are matrices of the form
\bdm
\g_{lk} = \left[
\begin{array}{ll}
0_{l,m+a_{21}+a_{31}-k} & \star \\
\star & \star
\end{array}
\right], \qquad (\x_2)_{ki} = \left[
\begin{array}{ll}
\star & \star \\
0_{k,1+a_{32}-i} & \star
\end{array}
\right].
\edm
\item[(ii)] $(\x_2,\x_3,\g)\nsim ((\x_2)_{ki},(\x_3)_i,\g_{lk})$, with $l/n > k\l_1+i\a_2$,
$0 \le k \le m+a_{21}+a_{31}-1$, $1 \le i \le 1+a_{32}$. Here $(\x_3)_i$ denotes a matrix
having zero entries on the last $i$ rows.
\item[(iii)] $(\x_2,\g) \nsim ((\x_2)_k,\g_{lk})$ with $l/n > k\l_1$, $0 \le k \le m+a_{21}+a_{31}-1$.
Here $(\x_2)_k$ denotes a matrix with zero entries on the last $k$ rows.
\end{enumerate}

\noindent \\
{\bf (8.1) Remark:} Let $W'$ be a matrix obtained by performing elementary row and
column operations on $W$. Assume that $W'$ has a zero submatrix with $a_{32}-1$
columns. Then the zero submatrix has at most $a_{21}$ rows.

\noindent \\
\emph{Proof:} We notice first that the $a_{21}$ nonzero entries of $W$ on each column
are linearly independent. The claim will follow if we show that the matrix $W_j$, obtained
by deleting the $j^{\text{th}}$ column of $W$ and those $i^{\text{th}}$ rows for which
$W_{ij}\neq 0$, has linearly independent rows. But, by construction, each row of $W_j$
is not zero and the nonzero entries of $W_j$ on each column are linearly independent.
This finishes the argument.

\noindent \\
As a direct consequence of the above remark we get the following:

\noindent \\
{\bf (8.2) Remark:} Assume that $\x_2 \sim (\x_2)_{k1}$. Then $k \le m + a_{21}$. \\

At the other extreme, we would like to know what is the largest $k$ for which
$\x_2 \sim (\x_2)_{k,a_{32}}$. For this we need the following analog of (8.1).
Its proof will be included in the proof of (8.5):

\noindent \\
{\bf (8.3) Remark:} Let $W'$ be as at (8.1). Then each column of $W'$ has $a_{21}$
linearly independent elements, in other words it spans $S^{d_1-d_2}V^*$.

\noindent \\
{\bf (8.4) Remark:} Assume that $\x_2 \sim (\x_2)_{k,a_{32}}$. Then $k \le m+a_{31}$. \\

Given integers $1\le i < j \le a_{32}$ we let $\o_{ij}$ be the number of nonzero rows
of the matrix made of the columns $i$ and $j$ of $W$. Let $\o$ be the smallest among the
numbers $\o_{ij}$. Note that $a_{21} < \o < 2a_{21}$. In fact, we have the obvious formula
\bdm
\o = 2 {d_1-d_2 +r \choose r} - {d_1-d_2 +r-1 \choose r}.
\edm

\noindent \\
{\bf (8.5) Remark:} Let $W'$ be as at (8.1). Assume that $W'$ has a zero submatrix
with two or more columns. Then the zero submatrix has at most $a_{31}-\o$ rows.

\noindent \\
\emph{Proof:} We have to show that any matrix made of two columns of $W'$
has $\o$ linearly independent rows. Let $W_j$ be the $j^{\text{th}}$ column of $W$.
Let $W_j'$ be a linear combination of the columns of $W$ of the form
\bdm
W_j' = W_j + \sum_{l < j} c_l W_l.
\edm
Given integers $1 \le p < q \le a_{32}$ we have to show that the matrix $W''$ made of the
columns $W_p'$ and $W_q'$ has $\o$ linearly independent rows.

We choose the lexicographic ordering on the monomials $X_0^{i_0} \cdots X_r^{i_r}$
that form a basis of $S^{d_2-d_3}V^*$.
We also choose the lexicographic ordering on the monomials giving a basis for $S^{d_1-d_3}V^*$.
We write $W$ relative to these orderings and we notice that, if the entry $W_{ij}$ of $W$ is
nonzero, then, for $i \le k$, $l \le j$, $(i,\, j)\neq (k,\, l)$, $W_{kl}$ is either zero or is larger
than $W_{ij}$ in the lexicographic ordering. This shows two things:
\begin{enumerate}
\item[(i)] if $W_{ij}\neq 0$, then $W_{ij}'$ is equal to $W_{ij}$ plus a linear combination of
monomials that are larger than $W_{ij}$ in the lexicographic ordering;
\item[(ii)] if $W_{ij}\neq 0$, then, for $i < k \le a_{31}$, $W_{kj}'$ is either zero or is a combination
of monomials larger than $W_{ij}$.
\end{enumerate}

\noindent
Performing on $W_j'$ row operations of the form $cR_k + R_i \lra R_i$, $i < k$, $c$ being
a scalar, we do not disturb properties (i) and (ii). Moreover, performing a certain sequence
of such operations, we can arrive at $W_j$. This proves remark (8.3).

To show that $W''$ has $\o$ linearly independent rows we proceed as follows. Performing, possibly, row operations on $W''$ of the kind mentioned above, we may assume that $W_p'=W_p$.
Now all we need to do is find $\o -a_{21}$ linearly independent elements among those $W_{iq}'$
for which $W_{ip}'=0$. But from (i) and (ii) we know that those $W_{iq}'$ for which $W_{iq}\neq 0$
are linearly independent. As there are at least $\o -a_{21}$ indices $i$ for which $W_{iq}\neq 0$
but $W_{ip}=0$, we are done.

\noindent \\
{\bf (8.6) Remark:} Assume that $\x_2 \sim (\x_2)_{ki}$ with $i \le a_{32}-1$. Then $k \le m+a_{21}
+ a_{31}-\o$. \\

Owing to the fact that $\x_2$ has $a_{21}+a_{31}$ linearly independent rows, we may assume that
$k \le m$ in (iii). Owing to the fact that $\x_3$ has $a_{32}$ linearly independent entries, we
may assume that $i=1$ in (ii) and, in view of (8.2), that $k \le m+a_{21}$.
According to (8.6), we may assume that $k \le m+a_{21}+a_{31}-\o$ in (i) if $2 \le i \le a_{32}-1$.
According to (8.4), we may assume that $k \le m+a_{31}$ in (i) if $i=a_{32}$.
Thus, in order to show that $(\x_2,\x_3,\g)$ is semistable, it is enough to show that
$\g \nsim \g_{lk}$
\begin{align*}
\text{with} \quad \frac{l}{n} & > k\l_1 && \text{and} && 0 \le k \le m, \\
\text{or with} \quad \frac{l}{n} & > k\l_1 + \a_2 && \text{and} && m < k \le m+a_{21}, \\
\text{or with} \quad \frac{l}{n} & > k\l_1 + i\a_2+ \a_3 && \text{and} &&
m+a_{21} < k \le m+a_{21}+a_{31}-\o, \\
& && && 2 \le i \le a_{32}-1,\\
\text{or with} \quad \frac{l}{n} & > k\l_1 + a_{32} \a_2 + \a_3 && \text{and} &&
m+a_{21}+a_{31}-\o < k \le m+a_{31}, \\
\text{or with} \quad \frac{l}{n} & > k\l_1 + (a_{32}+1)\a_2 +\a_3 && \text{and} &&
m+a_{31} < k \le m+a_{21}+a_{31} -1.
\end{align*}

\noindent \\
{\bf (8.7) Claim:} \emph{Let $W$ be the space of morphisms of sheaves on $\P^r$ of the form}
\bdm
\f : m\O(-d_1) \oplus \O(-d_2) \oplus \O(-d_3) \lra n\O.
\edm
\emph{Assume that $m < \o -a_{21}$, in other words assume that}
\bdm
m < {d_1-d_2+r-1 \choose r-1}.
\edm
\emph{Then for any nonsingular polarization $\L=(\l_1,\l_2,\l_3,\m_1)$ satisfying}
\begin{align*}
a_{21}  \l_1 <  \l_2 & < \frac{1-m\l_1 -a_{31}\l_1 + a_{32}a_{21}\l_1}{1+a_{32}}, \\
\l_2 & < 1-m\l_1, \\
\l_1 & \le \frac{1}{m+a_{21}} - \frac{1}{mn+a_{21}n} {d_3+r \choose r},
\end{align*}
\emph{the set of semistable morphisms admits a geometric quotient $W^{ss}(G,\L)/G$,
which is a quasiprojective variety}.

\noindent \\
\emph{Proof:} Let $\f$ be in $W^{ss}(G,\L)$. According to (2.6), we need to show that
$(\x_2,\x_3,\g(f))$ is semistable. We argue by contradiction. Assume that $\g(\f) \sim
\g_{lk}$ with $l/n > k\l_1$ and $0 \le k \le m$. Let $\psi = (\psi^1,\psi^2,\psi^3)$ denote
the truncated matrix consisting of the first $l$ rows of $\f$.
By assumption $\psi$ has kernel inside $k^m \oplus S^{d_1-d_2}V^* \oplus S^{d_1-d_3}V^*$ of dimension at least $m+a_{21}+a_{31}-k$, which is greater than $m+a_{21}$
and $m+a_{31}$ because, by hypothesis, $m< a_{21}$ and $m< a_{31}$.
This shows that $\Ker(\psi)$ intersects $S^{d_1-d_2}V^*$ and $S^{d_1-d_3}V^*$
nontrivially, forcing $\psi^2=0$ and $\psi^3=0$.
Moreover, replacing possibly $\f$ with an equivalent morphism, we may assume that there are at least $m-k$ linearly independent elements in the kernel of $\psi^1$,
viewed as a subspace of $k^m$. We get $\f \sim \f_{k,0,0}$, contradicting the semistability of $\f$.

Assume now that $\g(\f) \sim \g_{lk}$ with $l/n > k\l_1+\a_2$ and $m < k \le m+a_{21}$.
Note that automatically $l \ge l_{m,0,0}$. This excludes those $\g_{lk}$ with
$k < a_{21}$ because, as we saw above, the condition $k < a_{21}$ forces
$\psi^1=0$, $\psi^2=0$, so it yields $\f \sim \f_{m,0,0}$,
which is a contradiction. Thus $a_{21} \le k \le a_{21}+m$. We have
\bdm
\text{ker}(\psi^3) \ge a_{31} -k \ge a_{31}-a_{21}-m > a_{31}-\o \ge 0,
\edm
hence $\psi^3=0$. We cannot have $\psi^2=0$ because this would lead to the contradiction $\f \sim \f_{m,0,0}$. Thus the elements from $\Ker(\psi)$ project onto
$m+a_{21}-k$ linearly independent elements in $k^m$.
As $l/n > (k-a_{21})\l_1 + \l_2$, we obtain the contradiction $\f \sim \f_{k-a_{21},1,0}$.

Assume that $\g(\f) \sim \g_{lk}$ with $l/n > k\l_1 + i\a_2 + \a_3$, $m + a_{21} < k
\le m+a_{21}+ a_{31}-\o$ and $2 \le i \le a_{32}-1$. Notice that automatically
$l \ge l_{m,1,0}$, so we cannot have $\psi^3=0$, because this would lead to the contradiction $\f \sim \f_{m,1,0}$. Thus $\Ker(\psi)$ intersects $S^{d_1-d_3}V^*$
trivially, forcing
\bdm
m+a_{21}+a_{31}-k \le m +a_{21}, \quad \text{so} \quad
a_{31}\le k, \quad \text{so} \quad
a_{31} \le m+a_{21}+a_{31}-\o.
\edm
This contradicts the hypothesis on $m$.

We next examine the case $\g(\f) \sim \g_{lk}$ with $l/n > k\l_1 +a_{32}\a_2+ \a_3=
(k-a_{31})\l_1 + \l_3$ and $m+a_{21}+a_{31}-\o < k \le m+a_{31}$.
As before, $\psi^3 \neq 0$. Thus
\bdm
m+a_{21} \ge \text{ker}(\psi^1,\psi^2) \ge m+a_{21}+a_{31}-k \ge a_{21} > \o- a_{21}>m.
\edm
This shows that $k \ge a_{31}$ and that $\Ker(\psi)$ intersects $S^{d_1-d_2}V^*$
nontrivially. Replacing possibly $\f$ with an equivalent morphism, we may assume that
$\psi^2=0$ and that $\text{ker}(\psi^1) \ge m+a_{21}+a_{31}-k-a_{21}=m+a_{31}-k$.
We arrive at the contradiction $\f \sim \f_{k-a_{31},0,1}$.

The last situation we need to examine is $\g(\f) \sim \g_{lk}$ with
\bdm
\frac{l}{n}> k\l_1+(a_{32}+1)\a_2 + \a_3 = (k-a_{21}-a_{31})\l_1 + \l_2+\l_3
\edm
and $m+a_{31} < k \le m+a_{21}+a_{31}-1$. Notice that automatically $l \ge l_{m,0,1}$,
so $\Ker(\psi)$ intersects $S^{d_1-d_2}V^*$ and $S^{d_1-d_3}V^*$ trivially,
otherwise we would get the contradictions $\f \sim \f_{m,1,0}$ or $\f \sim \f_{m,0,1}$.

Assume that $(f,g)$ is a nonzero vector in $\Ker(\psi^2,\psi^3)$ viewed as a subspace
of $S^{d_1-d_2}V^* \oplus S^{d_1-d_3}V^*$. We saw above that both $f$ and $g$ must
be nonzero, so we can write $f=hf_1$, $g=hg_1$ with $f_1$, $g_1$ relatively prime
and
\bdm
\text{max} \{ 0,\, d_1-d_2-d_3 \} \le d = \text{deg}(h) \le d_1-d_2.
\edm
The rows of $(\psi^2,\psi^3)$ are of the form $(-g_1u,f_1u)$, hence they are vectors
in a vector space of dimension equal to $\text{dim}(S^{d_2+d_3+d-d_1}V^*)$.
Had we had the inequality $l \ge l_{m,0,0}+ \text{dim}(S^{d_2+d_3+d-d_1}V^*)$,
we would arrive at the contradiction $\f \sim \f_{m,0,0}$.
But this inequality follows from the condition
\bdm
n(a_{31}\l_1 +(a_{32}+1)\a_2 +\a_3) \ge \text{dim}(S^{d_3}V^*)
\edm
which is equivalent to the last condition from the statement of the claim.

In conclusion, $\Ker(\psi)$ intersects trivially $S^{d_1-d_2}V^* \oplus S^{d_1-d_3}V^*$,
hence its elements project onto linearly independent elements in $k^m$.
We get $m  \ge \text{ker}(\psi) \ge m+a_{21}+a_{31}-k$ forcing $k \ge a_{21}+a_{31}$.
We arrive at the contradiction $\f \sim \f_{k-a_{21}-a_{31},1,1}$.
This finishes the proof of the claim. The remaining conditions from the statement of the claim are there to ensure that $\l_3 > 0$ and $\a_3 > 0$. \\

According to corollary 7.2.2 in \cite{drezet-trautmann}, if $\a_2>0$, $\a_3>0$ and if
$\l_2 \ge a_{21}c_1(1,1)/n$, then the conclusion of (8.7) holds. The constant $c_1(1,1)$
can be computed as in the proof of lemma 9.1.2 (loc. cit.) One has
\bdm
c_1(1,1)= \frac{\text{dim}(S^{d_3}V)}{\text{dim}(S^{d_1-d_3}V)}.
\edm
Thus, according to Dr\'ezet and Trautmann, the conclusion of (8.7) holds under
the hypotheses
\begin{align*}
a_{21}\l_1 < \l_2 & < \frac{1-m\l_1-a_{31}\l_1 + a_{32}a_{21}\l_1}{1+a_{32}}, \\
\l_2 & < 1-m\l_1, \\
\l_2 & \ge \frac{a_{21}}{na_{31}} {r+d_3 \choose r}.
\end{align*}
Our result is not contained in Dr\'ezet and Trautmann's result.


\end{document}